\documentclass[10pt,leqno]{amsart}
\usepackage{amssymb}
\usepackage{mathrsfs}
\usepackage{color}
\usepackage{soul}

\newcommand\N{{\mathbb N}}
\newcommand\R{{\mathbb R}}

\def\AA{{\mathcal A}}
\def\BB{{\mathcal B}}
\def\CC{{\mathcal C}}
\def\DD{{\mathcal D}}
\def\EE{{\mathcal E}}
\def\FF{{\mathcal F}}

\def\HH{{\mathcal H}}

\def\KK{{\mathcal K}}

\def\QQ{{\mathcal Q}}
\def\RR{{\mathcal R}}

\def\TT{{\mathcal T}}

\def\VV{{\mathcal V}}

\def\ZZ{{\mathcal Z}}

\def\AAA{{\mathscr{A}}}
\def\BBB{{\mathscr{B}}}

\def\MMM{{\mathscr{M}}}

\def\eps{{\varepsilon}}

\newcommand{\wto}{\rightharpoonup}

\newtheorem{theo}{Theorem}

\newtheorem{lem}[theo]{Lemma}

\newtheorem{rem}[theo]{Remark}

\newtheorem{defin}[theo]{Definition}

\newcommand{\beqn}{\begin{equation}}
\newcommand{\eeqn}{\end{equation}}
\newcommand{\bear}{\begin{eqnarray}}
\newcommand{\eear}{\end{eqnarray}}
\newcommand{\bean}{\begin{eqnarray*}}
\newcommand{\eean}{\end{eqnarray*}}

\def\Nt{|\hskip-0.04cm|\hskip-0.04cm|}

\newcommand{\Black}{\color{black}}

\def\signge{\bigskip \begin{center} {\sc  Giani Ega\~na Fern\'andez \par\vspace{3mm}
Universidad de La Habana \par
Facultad de Matem\'atica y Computaci\'on\par
San L\'azaro y L, Vedado
CP 10400 C. Habana 
\par
CUBA\par\vspace{3mm}
e-mail:} \tt{gegana@matcom.uh.cu} \end{center}}

\def\signsm{\bigskip \begin{center} {\sc St\'ephane Mischler\par\vspace{3mm}
Universit\'e Paris-Dauphine \& IUF \par
CEREMADE, UMR CNRS 7534\par
Place du Mar\'echal de Lattre de Tassigny
75775 Paris Cedex 16\par
FRANCE\par\vspace{3mm}
e-mail:} \tt{mischler@ceremade.dauphine.fr} \end{center}}

\begin{document}

\title[Keller-Segel equation]
{Uniqueness and long time asymptotic for the Keller-Segel equation: the parabolic-elliptic case}

\author{G. Ega\~na, S. Mischler}

\begin{abstract}
The present paper deals with the parabolic-elliptic Keller-Segel  equation in the plane in
the general framework  of  weak (or ``free energy") solutions associated to initial datum with finite mass $M$, finite second moment and 
finite entropy. The aim of the paper is threefold:

(1) We prove the uniqueness of the ``free energy" solution on the maximal interval of existence $[0,T^*)$ with $T^*=\infty$
in the case when $M\le8\pi$ and $T^*< \infty$ in the case when $M> 8\pi$. The proof uses a DiPerna-Lions renormalizing argument which makes possible 
to get the  ``optimal regularity" as well as an estimate of the difference of two possible solutions in the critical $L^{4/3}$ Lebesgue norm 
similarly as for the $2d$  vorticity Navier-Stokes equation. 

(2) We prove immediate smoothing effect and, in the case $M < 8\pi$, we prove Sobolev norm bound uniformly in time for the rescaled solution (corresponding to the 
 self-similar variables).

(3) In the case $M < 8\pi$, we also prove weighted $L^{4/3}$ linearized stability of the self-similar profile and then universal optimal rate of convergence of the solution to the
self-similar profile. The proof is mainly based on an argument of enlargement of the functional space for  semigroup spectral gap.

\end{abstract}

\maketitle



 \smallskip


%
%
%
%
%
	
\vspace{0.3cm}



\bigskip


\tableofcontents



\section{Introduction} 
\label{sec:intro}
\setcounter{equation}{0}
\setcounter{theo}{0}


The aim of the paper is to prove uniqueness of  weak ``free energy" solutions to the 
the so-called parabolic-elliptic Keller-Segel  equation in the plane associated to initial datum with finite mass $M  \ge 0$, finite polynomial moment and 
finite entropy, and in the subcritical case $M <  8\pi$,  to prove optimal rate of convergence to self-similarity of these solutions.
 In  \cite{CM} our analysis will be extended to the  parabolic-parabolic Keller-Segel  equation in a similar context. 

\smallskip

The Keller-Segel (KS) system for chemotaxis describes the collective motion of cells that are attracted by a chemical substance that 
they are able to emit (\cite{P,KS}). We refer to \cite{MR2226917} and the references quoted therein for biological motivation and mathematical
introduction. In this paper we are concerned with the parabolic-elliptic KS model in the plane which takes
the form
\bear \label{eq:KS}
\partial_t f &=& \Delta f - \nabla (f \, \nabla c) \quad \hbox{in} \quad (0,\infty) \times \R^2, 
\\ \nonumber
c & := & - \bar\kappa  \,\,  = \,  - \kappa * f  \quad \hbox{in} \quad (0,\infty) \times \R^2,  
\eear
with $\kappa := {1 \over 2\pi}  \log |z| $,  so that in particular 
$$ 
- \nabla c = \bar\KK :=  \KK * f, \quad \KK := \nabla\kappa = {1 \over 2\pi} \, {z \over |z|^2}.
$$
Here $t \ge 0$ is  the time variable, $x \in \R^2$ is the space variable, $f=f(t,x) \ge 0$ stands for the {\it mass density of cells}  while $c = c(t,x) \in \R$ is the
{\it  chemo-attractant concentration} which solves the (elliptic) Poisson equation 
$-\Delta c =  f $ in $(0,\infty) \times \R^2$.

\smallskip

The evolution equation \eqref{eq:KS} is complemented with an initial condition 
\beqn\label{eq:KSt=0}
f(0,.) = f_0 \quad \hbox{in}\quad \R^2,
\eeqn 
where throughout this paper, we shall assume that 
\beqn\label{eq:BdsInitialDatum}
0 \le f_0 \in L^1_{2}(\R^2) , \quad f_0 \log f_0 \in L^1(\R^2).
\eeqn
Here and below for any weight function $\varpi : \R^2 \to \R_+$ we define the weighted Lebesgue space $L^p(\varpi)$ for $1 \le p \le \infty$ by 
$$
L^p(\varpi) := \{ f \in L^1_{loc}(\R^2); \,\, \| f \|_{L^p(\varpi)} := \| f \, \varpi \|_{L^p} < \infty \},
$$
as well as $L^1_+(\R^2)$ the cone of nonnegative functions of $L^1(\R^2)$. We also use the shorthand $L^p_k$, $k \ge 0$,  for the weighted  Lebesgue space associated to the polynomial growth weight function $\varpi(x) := \langle x \rangle^k$, $\langle x \rangle := (1 + |x|^2)^{1/2}$.  
 
\smallskip
The fundamental identities are that any solution to the Keller-Segel equation  \eqref{eq:KS}
satisfies at least formally the conservation of mass
\beqn\label{eq:MassConserv}
M(t) := \int_{\R^2} f(t,x) \, dx = \int_{\R^2} f_0(x) \, dx = : M,
\eeqn
the second moment equation 
\beqn\label{eq:2ndMoment}
M_2(t) := \int_{\R^2} f(t,x) \, |x|^2 \, dx = C_1(M) \, t + M_{2,0}, \quad M_{2,0} := \int_{\R^2} f_0(x) \, |x|^2 \, dx, 
\eeqn
$C_1(M):= 4 M \, \bigl( 1 - {M \over 8\pi} \bigr)$, 
and the free energy-dissipation of the free energy identity 
\beqn\label{eq:FreeEnergy}
\FF(t) + \int_0^t \DD_\FF(s) \, ds  = \FF_0,
\eeqn
where the free energy $\FF(t) = \FF(f(t))$, $\FF_0 = \FF(f_0)$ is defined by 
$$
\FF = \FF(f) := \int_{\R^2} f \log f dx + {1 \over 2} \int_{\R^2} f \, \bar\kappa  \, dx, 
$$
and the dissipation of free energy is defined by 
$$
\DD_\FF = \DD_\FF(f) := \int_{\R^2} f \, | \nabla (\log f) + \nabla \bar\kappa|^2 \, dx.
$$

It is worth emphasizing that the critical mass $M_* := 8\pi$ is a threshold because one sees from \eqref{eq:2ndMoment}
that there does not exist nonnegative and mass preserving solution when $M > 8\pi$ (the identity  \eqref{eq:2ndMoment} would imply that 
the second moment becomes negative in a finite time shorter than $T^{**} := 2\pi M_{2,0} / [M (8\pi-M)]$ which is in contradiction with the positivity of 
the solution).  

On the one hand, in the subcritical case $M < 8\pi$, 
thanks to the logarithmic Hardy-Littlewood Sobolev inequality (see e.g. \cite{MR1230930,MR1143664})
\beqn\label{eq:logHLSineq}
\forall \, f \ge 0, \quad 
\int_{\R^2} f (x) \log f(x) \, dx + {2 \over M} \int\!\! \int_{\R^2 \times \R^2} f(x) \, f(y) \, \log |x-y| \, dxdy   \ge C_2(M),
\eeqn
with $C_2(M) := M \, (1 + \log \pi - \log M)$, one can easily check   (see \cite[Lemma 7]{MR2226917}) that  
\beqn\label{eq:H<F}
\HH := \HH(f) = \int_{\R^2} f \log f \, dx \le C_3(M) \, \FF + C_4(M),
\eeqn
with $C_3(M) := 1/\bigl( 1 - {M \over 8\pi} \bigr) $, $C_4 (M) := C_3 (M) \, C_2(M) \, M/ (8\pi)$.
Then from \eqref{eq:H<F}  and  the very classical functional inequality (see for instance \cite[Lemma 8]{MR2226917}) 
\beqn\label{eq:H+<H}
\HH^+ := \HH^+(f) = \int_{\R^2} f (\log f)_+ dx \le \HH + { \frac 14} M_2 + C_5(M), 
\eeqn
with $C_5(M) := 2M\log (2\pi) + 2/e$, one concludes that \eqref{eq:MassConserv},  \eqref{eq:2ndMoment} and  \eqref{eq:FreeEnergy} 
 provide a convenient family of a priori 
estimates in order to define weak solutions. More precisely, we get
\bear\label{eq:3Identities}
&&\HH^+(f(t)) + M_2(f(t)) + C_3(M) \int_0^t \DD_\FF(f(s)) \, ds  \le  
\\ \nonumber
&&\le  C_3(M) \, \FF_0 + \frac 54 M_{2,0}  + 2 C_1(M) t +  C_4(M) + C_5(M), 
\eear
where the RHS term is finite under assumption \eqref{eq:BdsInitialDatum} on $f_0$, since 
\bear\label{eq:F<H} 
\FF_0 &\le&
\HH_0 + {1 \over 4\pi}  \int\!\! \int  f_0(x) \, f_0(y) \, (\log |x-y|)_+ \, dxdy 
\\ \nonumber
&\le&
\HH_0 + {1 \over 4\pi}  \int\!\! \int  f_0(x) \, f_0(y) \,   |x-y|^2 \, dxdy  \le \HH_0 + {1 \over \pi}  \, M \, M_{2,0} , 
\eear
with $\HH_0 := \HH(f_0)$. In other words, we have 
\beqn\label{eq:H+bdd}
\AAA_T(f) := \sup_{t \in [0,T]} \bigl\{  \HH^+(f(t)) + M_2(f(t)) \} + \int_0^T \DD_\FF(f(s)) \, ds  \le  C(T) \, \quad \forall \, T \in (0,T^*)
\eeqn
with $T^* = +\infty$ and a constant $C(T)$ which depends on $M$, $M_{2,0}$, $\HH_0$ and the final time $T$.

\smallskip
On the other hand, in the critical case $M = 8\pi$ and the supercritical case $M > 8\pi$, the above argument using
the logarithmic Hardy-Littlewood Sobolev inequality \eqref{eq:logHLSineq} fails, but one can however prove that 
\eqref{eq:H+bdd} holds with $T^*=+\infty$ when $M=8\pi$ and  that 
\eqref{eq:H+bdd} holds with some finial time $T^* \in (0,T^{**}]$ when $M > 8\pi$ (see \cite{BCM} for details as well as Remark~\ref{rem:Msupercritic} below).

\begin{defin}\label{def:sol} For any initial datum $f_0$ satisfying \eqref{eq:BdsInitialDatum} and any final time $T^* > 0$, we say that 
\beqn\label{eq:bddL1}
0 \le f \in L^\infty(0,T; L^1(\R^2)) \cap C([0,T); \DD'(\R^2)), \quad \forall \, T \in (0,T^*), 
\eeqn
is a weak solution to the Keller-Segel equation in the time interval $(0,T^*)$ associated to the initial condition $f_0$ whenever $f$ satisfies 
\eqref{eq:MassConserv},  \eqref{eq:2ndMoment} and 
\beqn\label{eq:FreeEnergyRelax}
\FF(t) + \int_0^t \DD_\FF(s) \, ds  \le \FF_0 \quad \forall \, t  \in (0,T^*),
\eeqn
as well as the Keller-Segel equation \eqref{eq:KS}-\eqref{eq:KSt=0} in the distributional sense, namely
\beqn\label{eq:KSweak}
 \int_{\R^2} f_0(x) \, \varphi(0,x) \, dx  = \int_0^{T^*} \int_{\R^2} f(t,x) \, \Bigl\{ (\nabla_x (\log f) +  \bar\KK )\cdot \nabla_x \varphi - \partial_t \varphi   \Bigr\} \, dxdt
\eeqn
for any  $\varphi \in C^2_c([0,T) \times \R^2)$. 
\end{defin}
 
 It is worth emphasizing that thanks to the Cauchy-Schwarz inequality, we have 
\bean
\int_{\R^2}  f \,  | \nabla_x (\log f) +  \bar\KK | \, dx 
\le M^{1/2} \, \DD_\FF^{1/2},
\eean
and the RHS of \eqref{eq:KSweak} is then well defined thanks to \eqref{eq:3Identities}. 

\medskip
This framework is well adapted for the 
existence theory.

\begin{theo} \label{theo:exist} For any initial datum $f_0$ satisfying \eqref{eq:BdsInitialDatum} there exists at least one weak solution on the 
time interval $(0,T^*)$ in the sense of Definition~\ref{def:sol} to the Keller-Segel equation  \eqref{eq:KS}-\eqref{eq:KSt=0} with $T^* = +\infty$ when 
$M \le 8\pi$ and $T^* < +\infty$ when $M >  8\pi$. 
  \end{theo}

We refer to \cite[Theorem 1]{MR2226917} for the subcritical case $M \in (0,8\pi)$ and to \cite{BCM} for the critical and supercritical cases $M \ge 8\pi$. 

\medskip
Our first main result establishes that this framework is also well adapted for the well-posedness issue. 

\begin{theo}\label{theo:uniq} For any initial datum $f_0$ satisfying \eqref{eq:BdsInitialDatum} there 
exists at most one weak solution in the sense of Definition~\ref{def:sol} to the Keller-Segel equation \eqref{eq:KS}-\eqref{eq:KSt=0}.
\end{theo}

 Theorem~\ref{theo:uniq} improves the uniqueness result proved in \cite{CLM} in the class of solutions $f \in C([0,T]; $ $ L^1_2(\R^2)) \cap L^\infty((0,T) \times \R^2)$
which can be built under the additional assumption $f_0 \in L^\infty(\R^2)$ (see also \cite{MR1654677} where a uniqueness result is established for a related model).  
Our proof follows a strategy introduced in \cite{FHM} for the 2D viscous vortex model. 
It is based on a DiPerna-Lions renormalization trick (see \cite{dPL}) which makes possible to get the optimal regularity of solutions for small time and 
then to follow the uniqueness argument introduced by Ben-Artzi for the 2D viscous vortex model (see \cite{MR1308857,MR1308858}). 
More precisely, we start proving the optimal regularity for short time $t^{1/4} \| f(t) \|_{L^{4/3}} \to 0$ as $t\to0$ for any weak solution $f$, and next we estimate the $L^{4/3}$-norm of the difference of two possible solutions  written in mild formulation. 
We emphasize that the $L^{4/3}$-space is critical for the Hardy-Littlewood-Sobolev inequality (see e.g. \cite[Theorem 4.3]{LL}) because it writes in that case
\beqn\label{eq:HLSineqCritic}
\bigl\|  f * \KK  \bigr\|_{L^{(4/3)'}\!(\R^2)} = \bigl\|  f * \KK  \bigr\|_{L^{4}(\R^2)}
\leq C   \, \| f \|_{L^{4/3}(\R^2)} , 
\eeqn
where $p' \in [1,\infty]$ is the conjugate exponent associated to $p \in [1,\infty]$  defined by $1/p + 1/p' = 1$. That last inequality is the key estimate in order to control the nonlinear term in \eqref{eq:KS}.  One probably could perform a similar argument with the $L^q$-norm, $q \ge 4/3$, see \cite{CM}.

\Black

\medskip
Next we consider the smoothness issue and the long time behaviour of solution for subcritical mass issue. 
For that last purpose it is convenient to work with self-similar variables. We  introduce the rescaled functions $g$ and $u$ defined by 
\beqn\label{eq:gt=}
g(t,x) := R(t)^{-2} f(\log R(t), R(t)^{-1} x), \quad 
u(t,x) :=  c(\log R(t), R(t)^{-1} x), 
\eeqn
with $R(t) := (1 + 2t)^{1/2}$. The rescaled parabolic-elliptic KS system reads 
\bear \label{eq:KSresc}
\partial_t g &=& \Delta g  + \nabla (g x - g \, \nabla u ) \quad \hbox{in} \quad (0,\infty) \times \R^2, 
\\ \nonumber
u &=& - \kappa *  g \quad \hbox{in} \quad (0,\infty) \times \R^2. 
\eear

Our second main result  concerns the regularity of the solutions. 

\begin{theo}\label{theo:ApostUnifSm} For any initial datum $f_0$ satisfying \eqref{eq:BdsInitialDatum}
the associated solution $f$ is smooth for positive time, namely $f \in C^\infty((0,T^* ) \times \R^2)$, 
and satisfies  the identity \eqref{eq:FreeEnergy} on $(0,T^*)$. 
Moreover, when $M < 8\pi$, the rescaled solution  $g$ defined by \eqref{eq:gt=}-\eqref{eq:KS}  satisfies the uniform in time moment estimate 
\beqn\label{eq:bddMkg}
\sup_{t \ge 0} M_k(g(t)) \le \max((k-1)^{k/2} M, M_k(f_0) )   \quad \forall \, k \ge 2,
\eeqn
with $M_k(g) := \| g \|_{L^1_k}$, as well as the  uniform in time regularity estimate for positive time 
\beqn\label{eq:bddW2infty}
\sup_{t \ge \eps} \| g (t,. ) \|_{W^{2,\infty}} \le \CC   \quad \forall \,\eps > 0,  
\eeqn
for some explicit constant $\CC$ which depends on $\eps$, $M$, $\FF_0$ and $M_{2,0} $. 
 \end{theo}
 
It is worth mentioning that $L^p$ bounds on $g$ for positive time and for $p \in [1,\infty)$  were known,  but non uniformly in time and as an a priori bound, while here  \eqref{eq:bddW2infty} is proved as an a posteriori estimate. 
Our proof is merely based on the same estimates as those established in  \cite{MR2226917},  on a bootstrap argument (using the DiPerna-Lions renormalization trick) and on the observation that 
 the rescaled free energy provides uniform in time estimates.  

\smallskip
From now on in this introduction, we definitively restrict ourself to the subcritical case $M < 8\pi$ and we focus on the long time asymptotic of the solutions. 
It has been proved in \cite[Theorem 1.2]{MR2226917}  that the solution given by Theorem~\ref{theo:exist} satisfies 
 \beqn\label{eq:gtoG}
 g(t,.) \to G \quad\hbox{in}\quad L^1(\R^2) \quad \hbox{as} \quad t \to\infty,
 \eeqn
 where $G$   is a solution to the rescaled stationary problem
  \bear \label{eq:statKSresc}
 &&  \Delta G  + \nabla (G x - G \, \nabla U ) = 0  \quad \hbox{in} \quad   \R^2, 
\\ \nonumber
&&  0 \le G, \quad \int_{\R^2} G \, dx = M, \quad U = - \KK * G . 
\eear
Moreover, the uniqueness of the solution $G$ to \eqref{eq:statKSresc} has been proved in \cite{MR2226917,BKLN}, see also \cite{MR2929020,MR2996772,CamposDolbeault2012}. From now on,  $G = G_M$  stands for the unique self-similar profile with same mass $M$ as $f_0$ 
and it is given in implicit form by 
 \beqn\label{eq:G=exp}
 G = M \, {e^{- G * \kappa  - |x|^2/2} \over \int_{\R^2} e^{- G * \kappa  - |x|^2/2} \, dx} , 
\eeqn
or equivalently $U = - \, G*\kappa$ satisfies
 \beqn\label{eq:U=exp}
 \Delta U + {M  \over \int_{\R^2} e^{U  - |x|^2/2} \, dx}   \, e^{U - |x|^2/2}  = 0.
\eeqn

\smallskip
Our third main result is about the convergence to self-similarity. 

\begin{theo}\label{theo:LongTimeRate} For any $M \in (0,8\pi)$, and any finite positive real numbers $\FF^*_0$, $k^* > 3$, $M^*_{k^*\!,0}$, 
there exists a  (non explicit) constant $C$  such that for any initial datum $f_0$ satisfying \eqref{eq:BdsInitialDatum} with 
$$
M_0(f_0) = M, \quad M_{k^*}(f_0) \le M^*_{k^*\!,0}, \quad \FF(f_0) \le \FF^*_0,
$$
the associated solution in self-similar variables $g$ defined by \eqref{eq:gt=}-\eqref{eq:KS}  satisfies  the optimal rate of convergence
$$
\| g(t,.) - G \|_{L^{4/3}} \le C \, e^{-t} \quad \forall \, t \ge 1,
$$
where $G$ stands for the self-similar profile with same mass $M$ as $f_0$.
\end{theo}

Let us emphasize that assuming only the second moment bound $M_2(f_0) < \infty$, 
the same proof leads to a not optimal rate of convergence to the self-similar profile, namely 
$\| g(t,.) - G \|_{L^{4/3}} \le C_\eta \, e^{-\eta t}$ for all $t \ge 1$ and for some $\eta \in (0,1)$, $C_\eta \in (0,\infty)$. It is likely that stronger 
moment assumption on the initial datum leads to the same optimal rate of convergence  in $L^q$-norm with larger values of $q$, but we do not follow that line of research in the present work. 

\smallskip
Theorem~\ref{theo:LongTimeRate} drastically improves some anterior results which establish the same exponential rate of convergence  for some particular class of initial data. 
On the one hand, for a radially symmetric initial datum with finite second moment it has been proved in \cite[Theorem 1.2]{MR2929020} the same convergence 
 in Wasserstein distance $W_2$ by a direct and nice entropy method. 
On the other hand, the same convergence in $L^2(e^{\nu |x|^2})$ norm, $\nu \in (0,1/4)$, has been recently proven to hold in \cite[Theorem 1]{CamposDolbeault2012} (see also \cite{MR2568716,MR2996772} for previous results in that direction) for any initial datum $f_0$ with mass $M \in (0,8\pi)$ and which satisfies (roughly speaking) the strong confinement condition $f_0 \le \tilde G$ for some self-similar profile $\tilde G$ associated to some mass $\tilde M \in [M,8\pi)$. 
In that last work~\cite{CamposDolbeault2012}, the uniform exponential stability (with optimal rate) of the linearized rescaled equation is established in $L^2(G^{-1/2})$ by the mean of the analysis of the associated linear operator in a well chosen (equivalent) Hilbert norm. The nonlinear exponential stability is then deduced from that linear stability together with an uniform in time estimate deduced from the strong confinement assumption made on the initial datum. 

\smallskip


Our proof follows a strategy of ``enlarging the functional space of semigroup spectral gap" initiated in \cite{Mcmp} for studying long time convergence to the equilibrium for the homogeneous Boltzmann equation,
and then developed in \cite{MMcmp,GMM,MR2832638,MR2821681,MM*} (see also \cite{MS*}) in the framework of kinetic equations and growth-fragmentation
equations. More precisely, taking advantage of the uniform exponential  stability of the linearized rescaled equation established in  \cite{CamposDolbeault2012}
in the small (strongly confining) space $L^2(G^{-1/2})$ (observe that $\log G (x) \sim -   |x|^2/2$ in the large position asymptotic) we prove that the same uniform exponential stability (with the same optimal rate) result holds in the larger space $L^{4/3}_\ell$, $\ell > 3/2$. It is worth emphasizing that the choice of the exponent $4/3$ is made in order to handle the singularity of the force field (thanks to the critical Hardy-Littlewood-Sobolev inequality \eqref{eq:HLSineqCritic}) while the choice of the moment exponent  $\ell > 3/2$ is made in order to have enough confinement and then to get the optimal rate.  We probably can perform a similar semigroup spectral gap analysis in a different space $L^q_\ell$, $q \ge 4/3$, $\ell \ge 3/2$,  but also probably at the cost of a stronger confinement (higher moment bound) assumption of the initial datum since at some point in the proof we use the embedding $L^q_\ell \subset L^\infty \cap L^1_{q\ell}$. Anyway, we do not follow that line of research in the present work. 
Next,  gathering the  long time convergence (without rate) to self-similarity \eqref{eq:gtoG}  with the estimates of Theorem~\ref{theo:ApostUnifSm},
 we obtain  that any solution reaches a small $L^{4/3}_\ell$-neighborhood of $G$ in finite time and we conclude to Theorem~\ref{theo:LongTimeRate}
 by nonlinear stability in $L^{4/3}_\ell$. It is worth emphasizing that it is only in that last nonlinear step that we use the a bit stronger initial (and then uniform in time) moment estimate \eqref{eq:bddMkg} with $k = k^* > 3$. 

\Black
 \medskip
 
 Let us end the introduction by describing the plan of the paper. 
 In Section~\ref{sec:aposteriori} we present some functional inequalities which will be useful in the sequel of the paper, 
we establish several a posteriori bounds satisfied by any weak solution, and we prove Theorem~\ref{theo:ApostUnifSm}.
Section~\ref{sec:uniqueness} is dedicated to the proof of the uniqueness result stated in Theorem~\ref{theo:uniq}. 
In Section~\ref{sec:LongTimeBehave} we prove the long time behaviour result as stated in Theorem~\ref{theo:LongTimeRate}. 

\medskip
\paragraph{\bf Acknowledgments} 

The authors gratefully acknowledge the support of the MADCOF ANR project (ANR-08-BLAN-0220).
E.G. would like to thank the CEREMADE at Universit\'e Paris-Dauphine for its kind hospitality 
in 2012 and 2013 where the work has been initiated and mostly written. 
He is also grateful to the MADCOF ANR project for the two several months grants that it provided to him. 
S.M. would like to thank the mathematics department of the Universidad de La Habana
for its hospitality  in summer 2013 where the present work has been concluded. 

%

 \section{A posteriori estimates - Proof of  Theorem~\ref{theo:ApostUnifSm} } 
\label{sec:aposteriori}
\setcounter{equation}{0}
\setcounter{theo}{0}

We start by presenting some elementary functional inequalities which will be of main importance in the sequel. 
The two first estimates are picked up from \cite[Lemma 3.2]{FHM} but are probably classical and the third one
is a variant of the Gagliardo-Nirenberg-Sobolev inequality.  

\begin{lem} \label{lem:FishInteg1}
For any $0 \le f \in L^1(\R^2)$ with finite mass $M$ and finite Fisher information
$$
I = I(f) := \int_{\R^2} {|\nabla f |^2 \over f}, 
$$
 there hold 
\bear
\label{eq:LpbdFisher} 
&\forall \; p\in [1,\infty), \quad \| f \|_{L^p(\R^2)} \leq C_p \,  M^{1/p} \,  I(f)^{1 - 1/p},\\
\label{eq:nablaLqbdFisher}
&\forall \; q\in [1,2), \quad
\| \nabla f \|_{L^q(\R^2)} \leq C_q  \, M^{1/q - 1/2} \,   I(f)^{ {3/2}-{1/q}}.
\eear
For any $0 \le f \in L^1(\R^2)$ with finite mass $M$, there holds 
\bear \label{eq:GNp+1}
&&\forall \; p \in [2,\infty) \quad
\| f \|_{L^{p+1}(\R^2)} \leq C_p  \, M^{1/(p+1)} \,   \| \nabla (f^{p/2}) \|_{L^2}^{2/(p+1)}.
 \eear

\end{lem}

For the sake of completeness we give the proof  below. 

\smallskip\noindent
{\sl Proof of Lemma~\ref{lem:FishInteg1}. } 
 We start with (\ref{eq:nablaLqbdFisher}).
Let $q \in [1,2)$ and use the H\"older inequality:
$$
\|\nabla f\|_{L^q}^q = \int \left|\frac{\nabla f}{\sqrt f}\right|^q f^{q/2}
\leq \left( \int \frac{|\nabla f|^2}{ f} \right)^{q/2} 
\left( \int  f^{q/(2-q)} \right)^{(2-q)/2} =  I(f)^{q/2} \, \|f \|_{L^{q/(2-q)}}^{q/2}. 
$$
Denoting by $q^\ast  =  2q/(2-q) \in [2,\infty)$ the Sobolev exponent associated to $q$ in dimension~$2$, we have,
thanks to a standard interpolation inequality and to the Sobolev inequality,
\bear\label{eq:Lq*bdFisher}
 \|f \|_{L^{q/(2-q)}} &=&  \|f \|_{L^{q^\ast / 2}} \le \| f \|_{L^1}^{1/ (q^\ast - 1)} \, 
\| f \|_{L^{q^\ast}}^{(q^\ast -2)/(q^\ast - 1)}
\\ \nonumber
& \le& C_q \, \| f \|_{L^1}^{1/(q^\ast - 1)} \, \| \nabla f \|_{L^{q}}^{(q^\ast -2)/(q^\ast - 1)}.
\eear
Gathering these two inequalities, it comes
$$
\|\nabla f\|_{L^q} \le  C_q \,  I(f)^{1/2} \, \| f \|_{L^1}^{1/ (2 (q^\ast - 1))} \, 
\| \nabla f \|_{L^{q}}^{(q^\ast -2)/(2(q^\ast - 1))},
$$
from which we   deduce \eqref{eq:nablaLqbdFisher}. 

\smallskip
We now establish \eqref{eq:LpbdFisher}.
For $p \in (1,\infty)$, write $p = q^\ast /2  =  q/(2-q)$ with  
$q:=2p/(1+p) \in [1,2)$ and use \eqref{eq:Lq*bdFisher} and \eqref{eq:nablaLqbdFisher} to get
\bean
\| f \|_{L^p} &\leq& 
 C_p \,  \| f \|_{L^1}^{ {1 \over q^\ast-1}+ {q^*-2 \over q^*-1} \, ( \frac1q-\frac12 ) }  
I(f)^{{q^*-2 \over q^*-1} \, (\frac32 -\frac1q)  },
\eean
from which one easily concludes. 

\smallskip
We  verify \eqref{eq:GNp+1}. From the Sobolev inequality and the Cauchy-Schwarz inequality, we have 
\bear\nonumber
\| w^{2(1+1/p)} \|_{L^1(\R^2)} 
&=&  \| w^{1+1/p} \|_{L^2(\R^2)}^2 \le \| \nabla (w^{1+1/p}) \|_{L^1(\R^2)}^2  \\
\label{eq:GNv-p+1}
&\le& (1+1/p)^2 \, \| w^{1/p} \|^2_{L^2} \,  \| \nabla w\|_{L^2(\R^2)}^2  
\eear
and we conclude to \eqref{eq:GNp+1} by taking $w := f^{p/2}$. 
\qed

\medskip
The proof of \eqref{eq:bddW2infty} in Theorem~\ref{theo:ApostUnifSm} is split into several steps that we present as some 
intermediate autonomous a posteriori bounds. 

\begin{lem}\label{lem:BdFisher} For any weak solution $f$ and any final time $T \in (0,T^*)$ there exists a constant $C := C(M,\AAA_T(f))$ such that 
\beqn\label{eq:boundIf}
\frac12 \int_0^T I(f(t)) \, dt \le C . 
\eeqn
In particular, in the subcritical case $M < 8\pi$ the constant $C$ only depends on $M$, $\HH_0$, $M_{2,0}$ and $T \in (0,\infty)$. 
\end{lem}

\smallskip\noindent
{\sl Proof of Lemma~\ref{lem:BdFisher}. } 
On the one hand, we write 
\bean
\DD_\FF(f) 
&=& \int f \, |\nabla (\log f + \bar\kappa)|^2
\\
&\ge& \int f \, |\nabla \log f |^2 + 2 \int \nabla f \cdot \nabla \bar\kappa 
=  I(f) - 2  \int f^2.
\eean
On the other hand, for any $A > 1$, using the Cauchy-Schwarz inequality and the inequality \eqref{eq:LpbdFisher}
 for
$p=3$, we have 
\bean
\int f^2 \, {\bf 1}_{f \ge A} 
&\le&
 \Bigl( \int f \, {\bf 1}_{f \ge A} \Bigr)^{1/2}  \Bigl( \int f^3\Bigr)^{1/2} 
\\ \nonumber
&\le&
 \Bigl( \int f \, {(\log f )_+\over \log A}  \Bigr)^{1/2}  \Bigl( C^3_3 \, M  \, I(f)^{2}\Bigr)^{1/2} ,
\eean
from what we deduce for $A = A(M, \HH^+(f))$  large enough, and more precisely taking $A$ such that $\log A = 16 \, \HH^+(f) \, C_3^3 \, M$, 
\beqn\label {eq:f2<I}
\int f^2 \, {\bf 1}_{f \ge A}  \le  C_3^{3/2} \, M^{1/2} \, {\HH^+(f)^{1/2} \over (\log A)^{1/2}} \, I(f) \le {1 \over 4} \, I(f).
\eeqn
Together with the first estimate, we find  
\bean
{1 \over 2} I(f) 
&\le& \DD_\FF(f)  + 2 \, \int f^2 \, {\bf 1}_{f \le A} 
\\
&\le& \DD_\FF(f)  + 2 \, M \, \exp (  16 \, \HH^+(f) \, C_3^3 \, M), 
\eean
and we conclude thanks to \eqref{eq:H+bdd} in the general case and thanks to  \eqref{eq:MassConserv}--\eqref{eq:F<H} in the subcritical case $M < 8\pi$. 
\qed  

\begin{rem}\label{rem:Msupercritic} As we have already mentioned we are not able to use the  logarithmic Hardy-Littlewood-Sobolev inequality \eqref{eq:logHLSineq} in the critical and supercritical cases. However, introducing the Maxwell function $\MMM := M \, (2\pi)^{-1} \exp(-|x|^2/2)$ of mass
$M$ and the relative entropy 
$$
H(h|\MMM) := \int_{\R^2} (h \, \log (h/\MMM) - h + \MMM) \, dx , 
$$
one classically shows that any solution $f$ to the Keller-Segel equation \eqref{eq:KS} formally satisfies 
\bean
{d \over dt} H(f(t)|\MMM) 
&=& 
- I(f(t)) + \int f(t)^2 + C_1/2 
\\
&\le& - I(f(t)) + MA + 
C_3^{3/2} \, M^{1/2} \, {\HH^+(f(t))^{1/2} \over (\log A)^{1/2}} \, I(f(t)) +C_1/2  \qquad (\forall \, A > 0)
\\
&=& - I(f(t)) + M  \exp \bigl( 4 C_3^3 \, M \, \HH^+(f(t)) )  +C_1/2 
\\
&=& - I(f(t)) + M  \exp \bigl\{ C_6 \, H(f(t)|\MMM)  \bigr\}  +C_1/2 ,
\eean
for a constant $C_6 = C_6(M)$ and where $C_1 = C_1(M)$ is defined in \eqref{eq:2ndMoment}. In the above estimates,  we have used \eqref{eq:f2<I}, we have made the choice $\log A := 4 C_3^3 \, M \, \HH^+(f(t))$ and
we have used a variant of inequality \eqref{eq:H+<H}. This differential inequality provides a local a priori estimate on the  relative entropy which is
the key estimate in order to prove local existence result for supercritical mass as well as global existence result for critical mass in 
\cite{BCM}. 
\end{rem}

As an immediate consequence of Lemmas~\ref{lem:FishInteg1} and \ref{lem:BdFisher}, we have 

\begin{lem}\label{lem:Bdf&K} For any $T \in (0,T^*)$, any  weak solution $f$ satisfies 
\bear \label{eq:wLqtLpx}
&& f \in L^{p/(p-1)}(0,T;L^p(\R^2)), \quad \forall \; p\in(1,\infty), 
\\
\label{eq:barKLqtLpx}
&&\bar\KK  \in L^{p/(p-1)}(0,T;L^{2p/(2-p)}(\R^2)),  \quad \forall \; p\in(1,2), 
\\
\label{bdd:dbbarK} 
&&\nabla_x \bar\KK   \in L^{p/(p-1)}(0,T;L^p(\R^2)) , \quad \forall \; p\in(2,\infty).
\eear

\end{lem}

\noindent
{\sl Proof of Lemma~\ref{lem:Bdf&K}. }  The bound \eqref{eq:wLqtLpx} is a direct consequence of \eqref{eq:boundIf} and \eqref{eq:LpbdFisher}. 
The bound \eqref{eq:barKLqtLpx} then follows from the definition of $K$, the Hardy-Littlewood-Sobolev inequality (see e.g. \cite[Theorem 4.3]{LL})
\beqn\label{eq:HLSineq}
\Bigl\|  {1 \over |z|} * f  \Bigr\|_{L^{2r/(2-r)}(\R^2)}
\leq C_r  \, \| f \|_{L^r(\R^2)}, \quad \forall \, r \in (1,2), 
\eeqn
with $r = p$ and  \eqref{eq:wLqtLpx}. Finally, from \eqref{eq:boundIf} and \eqref{eq:nablaLqbdFisher} we have 
$$
\nabla f \in L^{2q \over 3q-2} (0,T; L^q(\R^2)), \quad \forall \, q \in (1,2). 
$$
Applying the Hardy-Littlewood-Sobolev inequality \eqref{eq:HLSineq} to $\nabla_x \bar\KK = \KK * (\nabla_x f)$ with $r=q$, we get
$$
\nabla_x \bar\KK \in L^{2q \over 3q-2} (0,T; L^{2q \over 2-q}(\R^2)), \quad \forall \, q \in (1,2),
$$
which is nothing but  \eqref{bdd:dbbarK}. 
 \qed
 
\begin{lem}\label{lem:Bdbeta}  Any  weak solution $f$ satisfies 
\bear\label{eq:betaf}
&&\int_{\R^2} \beta(f_{t_1})  \, dx  + \int_{t_0}^{t_1} \int_{\R^2} \beta''(f_s) \, |\nabla f_s|^2 \, dx ds 
  \\ \nonumber
  &&\quad\le \int_{\R^2} \beta(f_{t_0})  \, dx  + \int_{t_0}^{t_1} \int_{\R^2} (\beta'(f_s) \, f^2_s - \beta(f_s) \, f_s)_+   \, dx ds ,
\eear
for  any times $0 \le t_0 \le t_1 < T^*$ and any renormalizing function $\beta : \R \to \R$ which is convex, piecewise of class $C^1$ and such that 
$$
|\beta (u)| \le C \, (1 + u \, (\log u)_+), \quad (\beta'(u) \, u^2 - \beta(u) \, u)_+  \le C \, (1 +  u^2) \quad \forall \, u \in \R.
$$
\end{lem}

\noindent
{\sl Proof of Lemma~\ref{lem:Bdbeta}. } We write 
$$
\partial_t f - \Delta_x f =   \bar\KK  \cdot \nabla_x f + f^2,
$$
and we split the proof into three steps. 
 
\smallskip\noindent
{\sl Step 1. Continuity. } 
Consider a mollifier sequence $(\rho_n)$ on $\R^2$, that is $\rho_n(x) := n^2 \rho(nx)$, $0 \le \rho \in \DD(\R^2)$, $\int \rho = 1$,  and introduce the
mollified function $f^n_t := f_t *_x \rho_n$. Clearly, $f^n  \in C([0,T), L^1(\R^2))$.
Using \eqref{eq:wLqtLpx} and  \eqref{bdd:dbbarK}, a variant of the commutation
Lemma \cite[Lemma II.1 and  Remark 4]{dPL} tells us that 
\beqn\label{eq:weps}
\partial_t f^n - \bar\KK \cdot \nabla_x f^n   - \Delta_x f^n = r^n, 
\eeqn
with
$$
r^n := (f^2) *\rho_n + (  \bar\KK \cdot \nabla_x f) * \rho_n -  \bar\KK \cdot \nabla_x
f^n \to f^2 \quad\hbox{in}\quad L^1(0,T; L^1_{loc}(\R^2)).
$$
The important point here is that $f^2, |\nabla_x \bar\KK| \,f \in L^1((0,T)
\times \R^2)$, thanks to~\eqref{bdd:dbbarK} and~\eqref{eq:wLqtLpx}.  

As a consequence,  the chain rule applied to the smooth function $f^n$ reads 
\beqn\label{eq:betaweps}
\partial_t \beta(f^n) =  \bar\KK \cdot \nabla_x \beta(f^n) + 
 \Delta_x \beta(f^n) -  \beta''(f^n) \, |\nabla_x f^n|^2 +
\beta'(f^n) \, r^n,
\eeqn
for any $\beta \in C^1(\R) \cap W^{2,\infty}_{loc}(\R)$ such that $\beta''$ is piecewise continuous and  
vanishes
outside of a compact set. Because the equation \eqref{eq:weps} with $\bar\KK$
fixed is linear, the difference 
$f^{n,k} := f^{n} - f^{k}$ satisfies \eqref{eq:weps} with $r^n$ replaced by $r^{n,k} := r^{n} - r^{k} \to 0$
in $L^1(0,T;L^1_{loc}(\R^2)$ and then also \eqref{eq:betaweps} (with again $f^n$ and $r^n$ changed in 
$f^{n,k}$ and $r^{n,k}$). In that last 
equation, we  choose $\beta(s) = \beta_1(s)$ where $\beta_A(s) = s^2/2$ for $|s| \le A$,  
$\beta_A(s) = A \, |s| - A^2/2$ for $|s| \ge A$ and we obtain 
for any non-negative function  $\chi \in C^2_c(\R^d)$,
\bean
&&\int_{\R^2} \beta_1(f^{n,k}(t,x)) \, \chi(x) \, dx  
\le 
\\
&&\le \int_{\R^2} \beta_1(f^{n,k}(0,x))\, \chi(x)  \, dx 
+  \int_{0}^{t} \!  \int_{\R^2}  |r^{n,k}(s,x)| \, 
\chi(x) \, dxds  
\\
&&
 + \int_{0}^{t} \!  \int_{\R^2}  \beta_1 (f^{n,k}(s,x)) \, \Bigl| - f \, \chi +   \Delta
\chi(x) - 
\bar\KK(s,x) \cdot \nabla \chi (x)\Bigr| \, dxds  , 
\eean
where we have used that $\hbox{div}_x \, \bar\KK = f$, that $|\beta_1'|
\le 1$  and that $\beta_1'' \ge 0$. In the last inequality, the RHS term converges to $0$ as $n,k$ tend to infinity. 
More precisely, $\beta_1(f^{n,k}(0) ) \to 0$ in $L^1(\R^2)$ because $f_0 \in L^1(\R^2)$;  $r^{n,k}  \to 0$ in  
$L^1(0,T; L^1_{loc}(\R^2))$ by the DiPerna-Lions commutation Lemma recalled above;  $ \beta_1(f^{n,k}) \bar\KK  \to 0$ in $L^1(0,T;L^1_{loc}(\R^2))$ 
because $\beta_1(s)\leq |s|$, because $f^{n,k}\to 0$ in $L^{3}(0,T,L^{3/2}(\R^2))$ 
by \eqref{eq:wLqtLpx} with $p=3/2$ and because $\bar\KK \in L^{6}(0,T;L^3(\R^2))$ by \eqref{eq:barKLqtLpx} with $p=6/5$; 
$ \beta_1(f^{n,k}) f  \to 0$ in $L^1(0,T;L^1 (\R^2))$ because again $\beta_1(s)\leq |s|$ and $f \in L^2((0,T) \times \R^2)$
by \eqref{eq:wLqtLpx} with $p=2$. 
All together, we get
$$
\sup_{t \in [0,T]} \int_{\R^2} \beta_1(f^{n,k}(t,x)) \, \chi(x) \, dx  
\, \mathop{\longrightarrow}_{n,k \to \infty} \, 0.
$$
Since $\chi$ is arbitrary, we deduce that there exists $\bar f \in C([0,\infty); L^1_{loc}(\R^2))$ so that 
$f^n \to \bar f$ in $C([0,T]; L^1_{loc}(\R^2))$, $\forall \, T > 0$. Together with the
convergence $f^n \to f$ in  $C([0,\infty);\DD'(\R^2))$ and the a priori bound \eqref{eq:3Identities}, we deduce  
that $f=\bar f$ and  
\beqn\label{eq:wcont}
f^n \to  f \quad \hbox{in} \quad C([0,T]; L^1(\R^2)), \quad \forall \, T > 0.
\eeqn

\smallskip\noindent
{\sl  Step 2. Linear estimates. } 
We come back to \eqref{eq:betaweps}, which implies, for all $0\le t_0<t_1$, all $\chi\in C^2_c(\R^2)$,
\bear\label{eq:betawnbis}
\int_{\R^2} \beta(f^{n}_{t_1}) \, \chi  \, dx  +    \int_{t_0}^{t_1} \! 
\int_{\R^2} 
\beta''(f^n_s) \, |\nabla_x f^n_s|^2 \, \chi\, dxds
=  \int_{\R^2} \beta( f^{n}_{t_0}) \, \chi  \, dx \quad
\\ \nonumber
+ \int_{t_0}^{t_1} \!  \int_{\R^2} \Bigl\{\beta'(f^n _s) \, r^n \, \chi +  \beta (f^n_s) \,   
 \Delta \chi  - 
\beta (f^n_s) \, \hbox{div}_x( \bar\KK  \chi) \Bigr\} \, dxds.  
\eear
Choosing $0 \le \chi \in C^2_c(\R^2)$ and 
$\beta \in C^1(\R) \cap 
W^{2,\infty}_{loc}(\R)$ such that $\beta''$ is non-negative and  vanishes
outside of a compact set,  and passing to the limit as $n\to\infty$,  we get 
\bear\label{eq:betafchi}
\int_{\R^2} \beta(f_{t_1}) \, \chi  \, dx  +    \int_{t_0}^{t_1} \! 
\int_{\R^2} \beta''(f_s) \, |\nabla_x f_s|^2 \, \chi\, dxds
\le  \int_{\R^2} \beta( f_{t_0}) \, \chi  \, dx \quad
\\ \nonumber
+ \int_{t_0}^{t_1} \!  \int_{\R^2} \Bigl\{ \Bigl[ \beta'(f) \, f^2   - \beta(f) \, f \Bigr] \, \chi 
+  \beta (f) \, \Bigl[  \Delta \chi  -  \bar\KK \cdot \nabla  \chi \Bigr] \Bigr\} \, dxds.  
\eear
By approximating $\chi\equiv 1$ by the sequence $(\chi_R)$ with $\chi_R(x) = \chi(x/R)$, $0 \le \chi \in \DD(\R^2)$, 
we see that the last term in \eqref{eq:betafchi} vanishes and we get \eqref{eq:betaf} in the limit $R \to \infty$ for 
any renormalizing function $\beta$ with linear growth at infinity.

\smallskip\noindent
{\sl  Step 3. superlinear estimates. } Finally, for any $\beta$ satisfying the growth condition as in the statement of the Lemma, 
we just approximate $\beta$ by an increasing sequence of smooth renormalizing functions $\beta_R$ with linear growth at infinity,
and pass to the limit in  \eqref{eq:betaf} in order to conclude.  \qed


\medskip
As a first consequence of Lemma~\ref{lem:Bdbeta}, we establish an estimate on the quantity 
\beqn\label{def:H2}
\HH_2(f) := \int_{\R^2} f \, (\widetilde \log f)^2 \, dx, \quad \widetilde \log \, u := {\bf 1}_{u \le e} + (\log u) {\bf 1}_{u > e}. 
\eeqn

\begin{lem}\label{lem:BdH2} 
For any weak solution $f$ and any time $T \in (0,T^*)$,  there exists a constant 
$C := C (M,T,\AAA_T)$ such that for any $0 \le  t_0 < t_1 \le T$ 
\beqn\label{eq:boundH2}
\HH_2(f(t_1)) \le \HH_2(f(t_0))+  C. 
\eeqn
\end{lem}

\noindent
{\sl Proof of Lemma~\ref{lem:BdH2}. }
We define the renormalizing function $\beta_K : \R_+ \to \R_+$, $K \ge e^2$,  by
$$
\beta_K(u) := u \, (\widetilde \log \,  u)^2 \,\,\,\hbox{if}\,\,\, u \le K, \,\,\,\quad
\beta_K(u) := (2+\log K) \, u \log u - 2 K \log K \, \,\,\hbox{if}\,\,\, u \ge K,
$$
so that  $\beta_K$ is convex and piecewise of class $C^1$, and moreover there holds  
$$
\beta_K'(u) u^2 - \beta_K(u) u \le 2 \, u^2 \, \widetilde \log \, u \, {\bf 1}_{u \le K} + 4 \log K \, u^2 \, {\bf 1}_{u > K}
$$
and
$$
\beta''_K(u) \ge 2 \, { \log  u \over u} \, {\bf 1}_{e \le u \le K} + (2 + \log K) \, {1 \over u} \, {\bf 1}_{u > K}.
$$
Defining now
$$
 \widetilde \log_K \, u := {\bf 1}_{u \le e} + (\log u) {\bf 1}_{e < u \le K} + (\log K) {\bf 1}_{u > K}, 
 $$
we deduce from \eqref{eq:betaf} that 
\bean
&& \int_{\R^2} \beta_K(f_{t_1}) \, dx  + \int_{t_0}^{t_1} \int_{\R^2} {  |\nabla f|^2  \over f}  \, (\widetilde \log_{K} f )\, {\bf 1}_{f \ge e} \, dx ds
  \\ \nonumber
  &&\quad\le  \int_{\R^2} \beta_K(f_{t_0}) \, dx +  4  \int_{t_0}^{t_1} \int_{\R^2}  f^2 \, \widetilde \log_K  f \,   dx ds .
\eean
Proceeding as in the proof of Lemma~\ref{lem:BdFisher}, we have for any $A \in (e,K)$
\bean
\int_{\R^2}  f^2 \, \widetilde \log_K  f  \, dx
&=&  \int_{\R^2}  f^2 \, \widetilde \log_K  f  \, {\bf 1}_{A \le K} \, dx
+   \int_{\R^2}  f^2 \, \widetilde \log_K  f  \, {\bf 1}_{A \ge K} \, dx
\\
&\le& (A \log A) \, M 
+  \Bigl( {\HH^+(f) \over \log A} \Bigr)^{1/2} \, \Bigl( \int_{\R^2}  (f \, \widetilde \log_K  f  )^3 \, dx \Bigr)^{1/2} ,
\eean
as well as thanks to inequality \eqref{eq:LpbdFisher} with $p=3$
\bean
\Bigl( \int_{\R^2}  (f \, \widetilde \log_K  f  )^3 \, dx \Bigr)^{1/2} 
&\le&
C_3^{3/2} \,  \Bigl( \int_{\R^2}  f \, \widetilde \log_K  f  \, dx \Bigr)^{1/2} 
 \int_{\R^2} { |\nabla  (f \, \widetilde \log_K  f ) |^2 \over f \, \widetilde \log_K  f}  \, dx  
 \\
&\le&
4 \, C_3^{3/2} \,  \Bigl( M + \HH^+(f) \Bigr)^{1/2}
\Bigl( \int_{\R^2}  {|\nabla f|^2 \over f} \, (\widetilde \log_K  f)  \, {\bf 1}_{ f \ge e} \, dx  + I(f) \Bigr). 
\eean
 The last three estimates together, we obtain for $A$ large enough and $K > A$
 \bean
&& \int_{\R^2} \beta_K(f_{t_1}) \, dx  \le  \int_{\R^2} \beta_K(f_{t_0}) \, dx +  4 \, T \,   (A \log A) \, M 
+  \int_{t_0}^{t_1} I(f_s)  ds,
\eean
from which \eqref{eq:boundH2} immediately follows by letting $K$ tends to $+\infty$ and using Lemma~\ref{lem:BdFisher}. 
\qed

\Black
\medskip
We now derive some $L^p$-norm estimate on the solutions to the KS equation. 

\begin{lem}\label{lem:BdLp} For any weak solution $f$, any time $T \in (0,T^*)$ and any $p \in[2,\infty)$ and $t_0 \in [0,T)$ such that $f_{t_0} \in L^p$,  there exists a constant 
$C := C (M,T,\AAA_T ,p, \| f_{t_0} \|_{L^p})$ such that 
\beqn\label{eq:boundLp}
\forall\, t_1 \in [t_0,T] , \qquad \| f (t_1) \|^p_{L^p}  + {1 \over 2} \int_{t_0}^{t_1} \| \nabla_x (f ^{p/2}) \|_{L^2}^2 \, dt \le C.
\eeqn
 \end{lem}

\noindent
{\sl Proof of Lemma~\ref{lem:BdLp}. }
 We define the renormalizing function $\beta_K : \R_+ \to \R_+$, $K \ge 2$,  by
$$
\beta_K(u) := \frac{u^p}p \,\,\, \hbox{if} \,\,\, u \le K, \quad 
\beta_K(u) := {K^{p-1} \over \log K} (u \, \log u - u) - {1 \over p'} K^p + {K^p \over \log K}    \,\,\, \hbox{if} \,\,\, u \ge K,
$$
so that  $\beta_K$ is convex and of class $C^1$, and moreover there holds  
$$
\beta_K'(u) u^2 - \beta_K(u) u \le {1 \over p'} u^{p+1} \, {\bf 1}_{u \le K} + 2 \,  K^{p-1}  \,  u^2 \, {\bf 1}_{u > K}, 
$$
as well as 
$$ 
\beta''_K(u)  = (p-1)\, u^{p-2}  \, {\bf 1}_{u \le K} + {K^{p-1} \over \log K}  \, {1 \over u} \, {\bf 1}_{u > K} .
$$

Thanks to Lemma~\ref{lem:Bdbeta}, we may write
\bean 
\int_{\R^2} \beta_K(f_{t_1})  \, dx  
+  {4 \over p' p}   \int_{t_0}^{t_1} \!  \int_{\R^2} |\nabla (f^{p/2})  |^2  \, {\bf 1}_{f \le K} \, dxds
+   {K^{p-1} \over \log K}   \int_{t_0}^{t_1} \!  \int_{\R^2} {|\nabla f  |^2 \over f}  \, {\bf 1}_{f \ge K} \, dxds
\\
\le  \int_{\R^2} \beta_K(f_{t_0})  \, dx    
+   {1\over p'}  \int_{t_0}^{t_1} \!  \int_{\R^2} f^{p+1} \,  {\bf 1}_{f \le K} \, dxds
  + 2 \, K^{p-1}    \int_{t_0}^{t_1} \!  \int_{\R^2}    f^2   \, {\bf 1}_{f \ge K}  \, dxds.
\eean

\Black


On the one hand, using the splitting $f = (f\wedge A) + (f-A)_+$, we have
\bean
\TT_1 
:=    \int_{\R^2} f^{p+1} \,  {\bf 1}_{f \le K} \, dx 
\le  2^p \, A^p \, M   +   2^p   \int_{\R^2} f_{A,K}^{p+1} \, dx ,
\eean
where we have defined $ f_{A,K} := \min( (f-A)_+, K-A)$, $K > A > 0$. Moreover, thanks to inequality \eqref{eq:GNp+1}
and the same trick as in the proof of Lemma~\ref{lem:BdFisher}, we have
\bean
\int_{\R^2} f_{A,K}^{p+1} \, dx 
&\le& C_p \,  \int_{\R^2} f_{A,K}  \, dx  \,  \int_{\R^2} |\nabla (f_{A,K}^{p/2}) |^2  \, dx
\\
&\le& C_p \,  {\HH^+(f) \over \log A}  \,  \int_{\R^2} |\nabla (f^{p/2}) |^2 \, {\bf 1}_{  f \le K}   \, dx. 
\eean
As a consequence, we obtain 
\bean
{1\over p'}  \int_{t_0}^{t_1} \TT_1  \, ds 
\le  {2^p \over p'} \, A^p \, M \, T  +  {1 \over p' p}  \int_{t_0}^{t_1} \! \int_{\R^2} |\nabla_x (f^{p/2}) |^2 \, {\bf 1}_{  f \le K}   \, dxds,
\eean
for $A = A(p,\AAA_T) > 1$ large enough. 

\smallskip
On the other hand,  thanks to the Sobolev inequality (line 2) and the Cauchy-Schwarz inequality (line 3), we have
\bean
\TT_2 
&:=& 2 K^{p-1}   \int_{\R^2}   f^2 \, {\bf 1}_{f \ge K} \, dx  
\le 8 \,  K^{p-1}    \int_{\R^2} (f-K/2)_+^2   \, dx  
\\
&\le&  8 \,  K^{p-1}  \Bigl( \int_{\R^2} |\nabla (f-K/2)_+|   \, dx  \Bigr)^2
= 8 \,  K^{p-1}   \Bigl( \int_{\R^2} |\nabla  f |  \, {\bf 1}_{f \ge K/2}  \, dx  \Bigr)^2
\\
&\le&   8 \,  K^{p-1}  \int_{\R^2} {|\nabla  f |^2 \over f}  \, {\bf 1}_{f \ge K/2}  \, dx \,  \int f \, {\bf 1}_{f \ge K/2} \, dx 
\\
&\le& 8 \,  K^{p-1}  \, \Bigl\{ {4 \over p^2} \int_{\R^2}    |\nabla  (f^{p/2}) |^2  \, \bigl( \frac 2K \bigr)^{p-1}  \, {\bf 1}_{f \le K}  + \int_{\R^2}   {|\nabla  f |^2\over f}  \, {\bf 1}_{  f \ge K}     \Bigr\}  {\HH_2(f) \over (\log(K/2))^2} . 
\eean
Recalling that from Lemma~\ref{lem:BdH2} we have 
\bean
\sup_{[t_0,T]} \HH_2(f) \le \HH_2(f_{t_0}) + C'  \le M + \| f_{t_0} \|_{L^p}^p  + C' =: C'',
\eean
we deduce 
\bean
 \int_{t_0}^{t_1} \TT_2  \, ds 
&\le&  {32 \, C'' \over (\log K)^2} \Bigl\{ {2^{p+1} \over p^2} \!
\int_{t_0}^{t_1}\!\! \! \int_{\R^2}  |\nabla  (f^{p/2}) |^2   \, {\bf 1}_{  f \le K}  \, dx ds
+ K^{p-1} \!\! \int_{t_0}^{t_1}\!\!  \! \int_{\R^2} {|\nabla  f |^2\over f}  \, {\bf 1}_{  f \ge K}   \, dxds \Bigr\}
\\
&\le&  {1 \over p'p}  \int_{t_0}^{t_1}\!\! \int_{\R^2}  |\nabla  (f^{p/2}) |^2   \, {\bf 1}_{  f \le K}  \, dxds +  {K^{p-1} \over \log K}  \int_{t_0}^{t_1}\!\! \int_{\R^2} {|\nabla  f |^2\over f}  \, {\bf 1}_{  f \ge K}   \, dxds,
\eean
for any $K \ge K^* = K^*(p,\AAA_T)> \max(A,4)$ large enough.  

\smallskip
All together, we have proved that for some constant $A$ and $K^*$ only depending on $p$, $T$, $\AAA_T$ and $f_{t_0}$,  and for any $K \ge K^*$ there holds 
 \bean 
\int_{\R^2} \beta_K(f_{t_1})  \, dx  
+  {2 \over p'}   \int_{t_0}^{t_1} \!  \int_{\R^2} |\nabla_x (f^{p/2})  |^2  \, {\bf 1}_{f \le K} \, dxds
\le  \int_{\R^2} \beta_K(f_{t_0})  \, dx  +    2^p \, A^p \, M \, T. 
\eean
We conclude to \eqref{eq:boundLp} by passing to the limit $K \to \infty$. 
\qed

\begin{lem}\label{lem:fsmooth} Any weak solution $f$ is smooth, that is  
$$
 f   \in C^\infty_b((\eps,T) \times \R^2), \,\, \forall \, \eps, T, \,\, 0< \eps < T < T^*, 
$$
so that in particular it is  a ``classical solution" for positive time.
\end{lem}

\noindent
{\sl Proof of Lemma~\ref{lem:fsmooth}. } For any time $t_0 \in (0,T)$ and any exponent $p \in (1,\infty)$,
there exists $t'_0 \in (0,t_0)$ such that $f(t'_0) \in L^p(\R^2)$ thanks to \eqref{eq:wLqtLpx},
from what we deduce using \eqref{eq:boundLp} on the time interval $(t'_0,T)$ that  
\beqn\label{eq:uniqvortexregL2nabla}
 f \in L^\infty(t_0,T; L^p(\R^2))
\quad \hbox{and}\quad 
\nabla_x f \in L^2((t_0,T) \times \R^2).
\eeqn
Next, by writing $\KK = \KK \, {\bf 1}_{|z| \le 1} + \KK \, {\bf 1}_{|z| \ge 1} \in L^{3/2}+L^\infty$, 
it is easily checked  $\| \KK * f \|_{L^\infty} \le C \, ( \| f \|_{L^3} + \| f \|_{L^1})$,
 and then   $\bar\KK \in L^\infty(t_0,T; L^\infty(\R^2))$ because of \eqref{eq:uniqvortexregL2nabla} and \eqref{eq:bddL1}.
We thus have
\beqn\label{eq:heatw}
\partial_t f + \Delta_x f = f^2 + \bar\KK \cdot \nabla_x f \in  L^2((t_0,T) \times \R^2),
\quad \forall \, t_0 >0, 
\eeqn
so that the maximal regularity of the heat equation in $L^2$-spaces (see Theorem X.11 stated in 
\cite{BrezisBook} and the quoted reference) provides the bound 
\beqn\label{eq:wMaximal}
f \in L^2(t_0,T; H^2(\R^2)) \cap L^\infty(t_0,T; H^1(\R^2)), \quad  \forall
\,t_0 >0.
\eeqn

Thanks to \eqref{eq:wMaximal}, an interpolation 
inequality and the Sobolev inequality, we deduce that
$\nabla_x f \in L^p ((t_0,T) \times \R^2)$ for any $1 < p < \infty$, whence $ \bar\KK \cdot \nabla_x
f \in L^p ((t_0,T) \times \R^2)$, for all $t_0 >0$. Then the maximal regularity
of the heat equation in $L^p$-spaces  
(see Theorem X.12 stated in \cite{BrezisBook} and the quoted references)  provides the bound 
\beqn\label{eq:wMaximalp}
\partial_t f, \nabla_x f  \in L^p((t_0,T) \times \R^2) , \quad \forall \, t_0 >0, 
\eeqn
and then the Morrey inequality implies the H\"olderian regularity $f \in C^{0,\alpha}((t_0,T) \times \R^2)$
for any $0< \alpha < 1$, and any $t_0 >0$.  Observing that the RHS term in \eqref{eq:heatw} has then also an H\"olderian regularity, we deduce that  
$$
\partial_t f, \partial_x f, \partial^2_{x_ix_j} f   \in C^{0,\alpha}_b((t_0,T) \times \R^2), \,\, \forall \, T, t_0; \,\, 
\,\,  0< t_0 < T < T^*, 
$$
thanks to the classical H\"olderian regularity result for the heat equation (see Theorem X.13 stated in \cite{BrezisBook} and the quoted references).
We conclude by (weakly) differentiating in time and space the equation \eqref{eq:heatw}, observing that the resulting RHS term is still a function with H\"olderian regularity, applying again \cite[Theorem X.13]{BrezisBook} and iterating the argument.  
\qed

\medskip\noindent
{\sl Proof of Theorem~\ref{theo:ApostUnifSm}. } 
We split the proof into  seven steps, many of them are independent from one another.

\smallskip\noindent
{\sl Step 1. } The regularity of $f$ has been yet established in Lemma~\ref{lem:fsmooth}.

\smallskip\noindent
{\sl Step 2. } First, we claim that the free energy functional $\FF$
is lsc in the sense that for any bounded sequence $(f_n)$ of nonnegative functions of $L^1_2(\R^2 )$ with same mass $M < 8\pi$ and such that $\FF(f_n) \le A  $
and $f_n \wto f$ in $\DD'(\R^2)$,  there holds 
\beqn\label{eq:Fsci}
0 \le f \in L^1_2(\R^2 ) 
\quad\hbox{and}\quad
\FF(f) \le \liminf \FF(f_{n}).
\eeqn
The proof of \eqref{eq:Fsci} is classical (see \cite{CLMP1,CLMP2,MR2226917}) and we just sketch it for the sake of completeness. 
Because of \eqref{eq:H<F} and  \eqref{eq:H+<H}, we have  $\HH^+(f_n) + M_2(f_n) \le A'$ for any $n \ge 1$, 
and we may apply the Dunford-Pettis lemma
which implies that $f_n \wto f$ in $L^1(\R^2)$ weak. Now, introducing the splitting $\FF = \FF_\eps + \RR_\eps$, $\FF_\eps = \HH + \VV_\eps$, with 
\bean
\VV_\eps(g) &:= & \frac12 \int \!\! \int_{\R^2 \times \R^2} g(x) \, g(y) \, \kappa(x-y) \, {\bf 1}_{|x-y| \ge \eps},
\\
\RR_\eps(g) &:=& \frac12 \int \!\! \int_{\R^2 \times \R^2} g(x) \, g(y) \, \kappa(x-y) \, {\bf 1}_{|x-y| \le \eps},
\eean
we clearly have that $\FF_\eps(f) \le \liminf \FF_\eps(f_{n})$ because $\HH$ is lsc and $\VV_\eps$ is continuous for the $L^1$ weak convergence.
On the other hand, using the convexity inequality $u v \le u \log u + e^v$ $\forall \, u > 0$, $v \in \R$ and the elementary inequality $(\log u)_- \le  u^{-1/2}$ 
$\forall \, u \in (0,1)$, we have for $\eps \in (0,1)$ and $\lambda > 1$ 
\bean
|\RR_\eps(g)|  
&=& {1 \over 4\pi} \int \!\! \int_{\R^2 \times \R^2} g(x) \, {\bf 1}_{g(x) \le \lambda} \, g(y) \, (\log |x-y|)_- \, {\bf 1}_{|x-y| \le \eps} 
\\
&& +  {1 \over 4\pi} \int \!\! \int_{\R^2 \times \R^2} g(x) \, {\bf 1}_{g(x) \ge \lambda} \, g(y) \, \log ( |x-y|^{-1} )  \, {\bf 1}_{|x-y| \le \eps} 
\\
&\le & {\lambda \over 4\pi}  \int_{\R^2} g(y) \, dy  \int_{|z|\le\eps} (\log |z|)_- \, dz
\\
&& +  {1 \over 4\pi}  \int_{\R^2} g(x) \, {\bf 1}_{g(x) \ge \lambda}
  \int \bigl\{ g(y) \log g(y) +   |x-y|^{-1} \bigr\} \, dy 
\\
&\le & {\lambda \over 3} \, M \, \eps^{3/2}
+  
 {1 \over 4\pi} \, {\HH^+(g) \over \log \lambda} \,  \bigl\{ \HH^+(g) +  2\pi\eps \bigr\} ,
\eean
and we get that $\sup_n |\RR_\eps(f_n)| \to 0$ as $\eps \to 0$ from which we conclude that $\FF$ is lsc. Now, we easily deduce that the free energy identity \eqref{eq:FreeEnergy} holds. Indeed, since $f$ is smooth for positive time, for any fixed $t \in (0,T^*)$ and any given sequence $(t_n)$ of positive real numbers which decreases to $0$, we clearly 
have 
$$
\FF(f(t_n)) = \FF(t) + \int_{t_n}^t \DD_\FF(f(s)) \, ds. 
$$  
Then, thanks to the Lebesgue convergence theorem,  the lsc property of $\FF$ and the fact that $f(t_n) \wto f_0$ weakly in $\DD'(\R^2)$, we deduce 
from the above free energy identity  for positive time that 
$$
\FF(f_0) \le \liminf_{n\to\infty} \FF(f(t_n)) \le \lim_{n\to\infty} \bigl\{ \FF(t) + \int_{t_n}^t \DD_\FF(f(s)) \, ds \bigr\} =  \FF(t) + \int_{0}^t \DD_\FF(f(s)) \, ds.
$$ 
Together with the reverse inequality \eqref{eq:FreeEnergyRelax} we conclude to \eqref{eq:FreeEnergy}. 

\smallskip\noindent
{\sl Step 3. } From now on, we assume that $M < 8\pi$ is subcritical and we prove the uniform in time estimates \eqref{eq:bddMkg}
and \eqref{eq:bddW2infty}. 
We start with the a priori additional moment estimate  \eqref{eq:bddMkg}. Because we will show the uniqueness of solution without using that 
additional moment estimates, these ones are rigorously justified thanks to a standard approximation argument, see \cite{MR2226917} for details. 
Denoting $g$ the rescaled solution \eqref{eq:gt=} and 
$$
 M_k := \int_{\R^2} g(x) \, |x|^k \, dx
$$
we compute with   $\Phi(x) = |x|^k$, $k \ge 2$,  thanks to the antisymmetry of the kernel and the H\"older inequality
\bean
{d \over dt} M_k 
&=& k^2 \, M_{k-2} - k \, M_k -{1 \over 2\pi} \int_{\R^2} \Phi'(x) g(t,x) \int_{\R^2} g(t,y) \, {x-y \over |x-y|^2} \, dy dx
\\
&=&  k^2 \, M_{k-2} - k \,  M_k 
\\
&&- {1 \over 4\pi} \int_{\R^2}\! \int_{\R^2} g(t,y) g(t,x) \, ( \Phi'(x) -  \Phi'(y)) {x-y \over |x-y|^2} \, dy dx
\\
&\le&k^2\, M^{2/k} \, M_k^{1 - 2/k} - k \,  M_k  ,
\eean
from which we easily conclude that \eqref{eq:bddMkg}
holds.  

\smallskip\noindent
{\sl Step 4. } Defining the rescaled free energy $\EE(g)$ and the associated dissipativity of rescaled free energy $\DD_\EE(g)$ by
\bear\label{def:EE}
\EE(g) &:=& \int g ( 1 +  \log g)  + {1 \over 2} \int g|x|^2 + { {1 \over 4 \pi} } \int\!\!\int g(x) g(y) \log |x-y| \, dxdy 
\\ \label{def:DDEE}
\DD_\EE(g) &:=&  \int  g \Bigl|\nabla \Bigl(\log g+ \frac{|x|^2}2 + \kappa * g \Bigr)\Bigr|^2,
\eear
we have that any solution $g$ to the rescaled equation \eqref{eq:KSresc} satisfies 
\beqn\label{eq:dissipRescE}
{d \over dt} \EE(g) + \DD_\EE(g) = 0 \quad \hbox{on} \quad [0,\infty). 
\eeqn
On the one hand, as for \eqref{eq:H<F}, the following functional inequality 
\beqn\label{eq:H<E}
\int g   \log g   + {1 \over 2} \int g|x|^2 \le C_3(M) \, \EE (g)+ C_4(M) \quad \forall \, g \in L^1_+(\R^2)
\eeqn
holds, and together with \eqref{eq:H+<H},  we find 
\beqn\label{eq:H+<E}
\int g (\log g)_+   + {1 \over 4} \int g|x|^2 \le C_3(M) \, \EE (g) + C_7(M)  \quad \forall \, g \in L^1_+(\R^2),
\eeqn
where $C_7 := C_4 + C_5$. As a consequence of \eqref{eq:dissipRescE} and \eqref{eq:H+<E}, we get the uniform in time 
upper bound on the rescaled free energy for the solution $g$ of  \eqref{eq:KSresc} 
\beqn
\sup_{t \ge 0} \int g_t (\log g_t)_+ + \frac14 \int g_t |x|^2 \le  C_3(M) \, \EE(f_0) + C_7(M).
\eeqn

\smallskip\noindent
{\sl Step 5. } As in the proof of Lemma~\ref{lem:BdLp}, we easily get that the rescaled solution $g$ of the rescaled equation \eqref{eq:KSresc} 
satisfies for any $p \in [2,\infty)$
\bean
{d \over dt} \| g \|_{L^p}^p 
+  {4 \over p'}  \| \nabla (g^{p/2}) \|^2_{L^2} 
&=& 2 \, ( p-1)   \,  \| g \|_{L^p} ^p  + (p-1)  \| g \|_{L^{p+1}} ^{p+1} 
\\
&\le& 2 \, ( p-1)   \, M   + 3(p-1)  \| g \|_{L^{p+1}} ^{p+1} . 
\eean
Writing $s = s \wedge A + (s-A)_+$, so that  $s^{p+1} \le 2^{p+1}(s \wedge A)^{p+1} + 2^{p+1}(s-A)_+^{p+1} $,
and using the Gagliardo-Nirenberg-Sobolev type inequality \eqref{eq:GNv-p+1} in order to get 
\bean
\int (g-A)_+^{p+1} &\le& C_p \, \int |\nabla (g-A)_+^{p/2}|^2 \, \int (g-A)_+ 
\\
&\le& C_p \, \int |\nabla (g^{p/2}) |^2 \, {\HH^+(g) \over \log A}
\eean
for any $A > 1$, we deduce  
\bean
{d \over dt} \| g \|_{L^p}^p + \| \nabla (g^{p/2}) \|^2_{L^2}  
&\le& 2pM + 3p2^{p+1} A^p \, M +  3p2^{p+1}   \int (g-A)_+^{p+1} \\
&\le& C_8(M,p,A)+ C_p \, {\HH_+(g) \over \log A}  \| \nabla (g^{p/2}) \|^2_{L^2}. 
\eean
Taking $A$ large enough,  we obtain 
\beqn\label{eq:gLpunif}
{d \over dt} \| g \|_{L^p}^p + {1 \over 2} \| \nabla (g^{p/2}) \|^2_{L^2}  
\le  C_{9}(M,p,\EE_0).
\eeqn
Using the Nash inequality 
$$
\| w \|_{L^2(\R^2)}^2 \le C_N \, \|w \|_{L^1(\R^2)} \, \| \nabla w \|_{L^2(\R^2)}
$$
with $w := g^{p/2}$, we conclude with 
\bean
{d \over dt} \| g \|_{L^p}^p + {1 \over C_N^2} \, \| g \|_{\, L^{p/2}}^{-p} \, \| g \|_{L^p}^{2p} 
\le  C_{9}(M,p,\EE_0).
\eean
Defining $u(t) := \| g(t) \|_{L^p}^{p}$ first with  $p=2$, so that  $\| g(t) \|_{L^{p/2}}^{p/2} = M$, we recognize  the classical 
nonlinear ordinary differential inequality 
$$
u' + c \, u^2 \le C \quad \hbox{on} \quad (0,\infty),
$$
for some constants $c$ and $C$ (which only depend on $M$ and $\EE_0$) from which we deduce the 
bound
\beqn
\label{eq:bddL2}
\forall \, \eps > 0 \,\, \exists \, \CC = \CC(\eps,c,C) \quad
\sup_{t \ge \eps} \| g(t) \|_{L^p}^{p} \le \CC,
\eeqn
with $p=2$.  
In order to get the same  uniform estimate \eqref{eq:bddL2} in all the Lebesgue spaces $L^p$, $p \in (2,\infty)$,  we may proceed by iterating the same argument as above with the choice $p = 2^k$, $k \in \N^*$. Coming back to \eqref{eq:gLpunif} with $p=2$, we also deduce that for any $\eps, T > 0$ there exists $\CC = \CC(\eps,T, \EE_0)$ so that 
$$
\sup_{t_0 \ge \eps} \int_{t_0}^{t_0+T} \| \nabla g (s) \|_{L^2(\R^2)}^2 \, ds \le \CC.
$$

\smallskip\noindent
{\sl Step 6. }  The function $g_i := \partial_{x_i} g$ satisfies 
$$
\partial_t g_i - \Delta g_i - \nabla (x g_i) = g_i + 2 \, g \, g_i - \partial_{x_i} ( \nabla u \cdot \nabla g),
$$
from which we deduce that 
\bear\label{eq:ODIgip}
&&{d \over dt} \int |g_i|^p + p(p-1) \int |\nabla g_i|^2 |g_i| ^{p-2} \le
\\ \nonumber
&&\le  (3p-2) \int |g_i|^p + 2p \int g \,  |g_i|^p  + p \int   \partial_{x_i} ( \nabla u \cdot \nabla g)   \, g_i \, |g_i|^{p-2}.
\eear
For $p=2$, we have for any $t \ge \eps$
\bean
\TT(t) &:=& 4 \int g \,  |g_i|^2  + 2 \int   \partial_{x_i} ( \nabla u \cdot \nabla g)   \, g_i  
\\
&\le& 4 \| g \|_{L^3} \, \| g_i \|_{L^3}^2  + 2 \,  \|  \nabla u \cdot \nabla g \|_{L^2}^2 + {1 \over 2} \,  \| \partial_i g_i \|_{L^2}^2   
\eean
thanks to the  H\"older inequality, an integration by part and the Young inequality. Next, 
we have for any $t \ge \eps$
\bean
\TT(t)
&\le& C_1  \, \| g_i \|_{L^2}^{4/3} \| \nabla g_i \|_{L^2}^{2/3}  + C_2  \|    \nabla g \|_{L^2}^2 + {1 \over 2} \,  \| \nabla  g_i \|_{L^2}^2   
\eean
where we have used  the classical Gagliardo-Nirenberg inequality (see (85) in \cite[Chapter IX]{BrezisBook} and the 
quoted references) 
\beqn\label{eq:GNineq}
\| w \|_{L^r(\R^2)} \le C_{GN} \, \| w \|_{L^q(\R^2)}^{1-a} \, \| \nabla w \|_{L^2(\R^2)}^a,
\quad a= 1- {q\over r}, \,\,\, 1 \le q \le r < \infty,
\eeqn
with $w:=g_i$, $r=3$, $q=2$, the uniform bound established in step 5 and the fact that 
 $\nabla u = - \KK * g  \in L^\infty((\eps,\infty) \times \R^2)$ thanks to the same argument as in the proof of Lemma~\ref{lem:fsmooth}. 
Last, by the Young inequality we get  for any $t \ge \eps$
\bean
\TT(t) &\le& {2 \over 3} \, C_1^{3/2}   \, \| g_i \|_{L^2}^2 + {1 \over 3} \,  \| \nabla g_i \|_{L^2}^2  + C_2  \|    \nabla g \|_{L^2}^2 + {1 \over 2} \,  \| \nabla  g_i \|_{L^2}^2,   
\eean
from which we deduce from \eqref{eq:ODIgip}
\bean
{d \over dt} \int |g_i|^2 +   \int |\nabla g_i|^2  \le C_3 \,  \|    \nabla g \|_{L^2}^2
  \quad\hbox{on}\quad (\eps,\infty),
\eean
with $C_3 :=   4 + {2 \over 3} \, C_1^{3/2}  +  C_2$. 
Remarking that for any fixed $\eps \in (0,1)$ and any $t_1 \ge 2\eps$,  we may define $t_0 \in (t_1-\eps,t_1)$ so that 
\bean
\| \nabla g (t_0) \|_{L^2}^2 
&=& \inf_{(t_1-\eps,t_1)} \| \nabla g \|_{L^2}^2 \le {2 \over \eps} \int_{t_1-\eps}^{t_1} \| \nabla g (s) \|_{L^2}^2 \, ds \le C_4
\eean
thanks to the bound established at the end of step 5, we deduce from the above differential inequality that 
$$
\| g_i(t_1) \|_{L^2}^2  \le\| g_i(t_0) \|_{L^2}^2  + C_3 \int_{t_0}^{t_1} \| \nabla  g (s) \|_{L^2}^2 \, ds  \le C_5,
$$
where again $C_5 := C_4 + C_3 \, C_4 \, \eps /2$ only depends on $\eps$, $M$ and $\EE_0$. 
Coming back to the above differential inequality again, we easily conclude that for any $\eps > 0$, there exists a constant $\CC_\eps = \CC(\eps,M,\EE_0)$
so that 
\beqn\label{eq:gH1H2unif}
\sup_{t \ge \eps} \Bigl\{ \| \nabla g (t) \|_{L^2}^2 + \int_{t}^{t+1}  \| D^2 g (s) \|_{L^2}^2 \Bigr\} \le \CC_\eps.
\eeqn 
 
 \smallskip\noindent
{\sl Step 7. } Starting from  the differential inequality \eqref{eq:ODIgip} for $p \in (2,\infty)$ and using 
 the Morrey-Sobolev inequalities
$$
\| g \|_{L^\infty} \le C \, \| g \|_{H^2} 
\quad\hbox{and}\quad
\| D^2 u \|_{L^\infty} \le C \,  \| D^2 u \|_{H^2} \le C \, \| g \|_{H^2},
$$
 we easily get 
\bean
\frac 1p  {d \over dt} \int |\nabla g|^p  
&\le&  C \, (1 + \| g \|_{L^\infty} +  \| D^2 u \|_{L^\infty} )  \int |\nabla g|^p  
\\
&\le&  C \, (1 + \| g \|_{H^2})  \int |\nabla g|^p   \quad\hbox{on}\quad (\eps,\infty),
\eean
from which we deduce for any $t_1 \ge t_0 \ge \eps$
$$
 \| \nabla g(t_1) \|_{L^p} \le  \| \nabla g(t_0) \|_{L^p} \, \exp \Bigl( \int_{t_0}^{t_1} C \, (1 + \| g(s) \|_{H^2}) \, ds \Bigr).
 $$
Now, arguing similarly as in step 6, we deduce from  the above time integral inequality, the Sobolev inequality 
$\| \nabla g \|_{L^p} \le C_p \, \| g \|_{H^2}$ for $p \in [2,\infty)$ and the already established bound \eqref{eq:gH1H2unif}, 
that for any $\eps > 0$, there exists a constant $\CC_\eps = \CC(\eps,M,\EE_0,p)$
so that 
\beqn\label{eq:gW1p}
\sup_{t \ge \eps}   \| \nabla g (t) \|_{L^p}  \le \CC_\eps.
\eeqn 

 \smallskip\noindent
{\sl Step 8. } Iterating twice the arguments we have presented in steps 6 and 7, it is not difficult to prove 
$$
\sup_{t \ge \eps} \| g (t,. ) \|_{W^{3,p}} \le \CC   \quad \forall \,\eps > 0, \,\, p \in [2,\infty), 
$$
for some constant $\CC = \CC(\eps,p,M,\FF_0,M_{2,0} )$ from which \eqref{eq:bddW2infty} immediately follows.
 \qed
 
 \section{Uniqueness - Proof of Theorem~\ref{theo:uniq} } 
\label{sec:uniqueness}
\setcounter{equation}{0}
\setcounter{theo}{0}

  We split the proof into two steps. We recall that from Theorem~\ref{theo:ApostUnifSm} we already know that
  $\| f \|_{L^2} \in C^1(0,T)$ and $\| f \|_{L^p} \in L^\infty(t_0,T)$ for any $0 < t_0 < T  < T^*$ and any $p \in [1,\infty]$.

\smallskip\noindent
{\sl Step 1. } We establish our new main estimate, namely that any weak solution satisfies
\beqn\label{eq:L43to0}
t^{1/4} \| f (t,.) \|_{L^{4/3}} \to 0 \quad\hbox{as}\ t \to 0.
\eeqn

First, from \eqref{eq:KS} and the regularity of the solution, we have 
$$
{d \over dt} \| f \|_{L^2}^2 + 2 \| \nabla_x f \|_{L^2}^2   =   \| f \|_{L^3}^3 \quad\hbox{on}\quad (0,T).
$$
As in the proof of Lemma~\ref{lem:BdLp}, we deduce that 
$$
{d \over dt} \| f \|_{L^2}^2 +{1 \over 2}  \| \nabla_x f \|_{L^2}^2  \le A^2 M\quad\hbox{on}\quad (0,T)
$$ 
for $A$ large enough. Thanks to the Nash inequality
$$
\| f \|_{L^2}^2 \le C \, M \, \| \nabla f \|_{L^2},
$$
we thus obtain  $$
{d \over dt} \| f \|_{L^2}^2 + c_M \| f \|_{L^2}^4 \le A^2 M \quad\hbox{on}\quad (0,T).
$$ 
It is a classical trick of ordinary differential inequality to deduce that there exists a constant $K$ (which only depends on $c_M$,  $A^2 M$ and $T$) 
so that 
\beqn\label{eq:L2tC}
t \, \|  f (t, .) \|_{L^2}^2 \le K \  \quad \forall \, t \in  (0,T).
\eeqn
We now prove \eqref{eq:L43to0} from \eqref{eq:L2tC} and an interpolation argument. 
On the one hand, introducing the notation $\widetilde{\log}_+ f := 2 + (\log f)_+$, we use the H\"older inequality in order to get 
\bean
\int f^{4/3} 
&=& \int f^{2/3} \, (\widetilde{\log}_+ f)^{2/3} \, f^{2/3} \,(\widetilde{\log}_+ f)^{-2/3} 
\\
&\le& \Bigl( \int f  \, \widetilde{\log}_+ f \Bigr)^{2/3} \, \Bigl( \int f^{2} \, (\widetilde{\log}_+ f)^{-2} \Bigr)^{1/3},
\eean
or in other words and using a similar estimate as \eqref{eq:H+<H}
\beqn\label{eq:L43bdd}
\| f \|_{L^{4/3}} \le C(\HH(f),M_2(f)) \, \Bigl( \int f^{2} \, (\widetilde{\log}_+ f)^{-2} \Bigr)^{1/4}.
\eeqn
On the other hand, we observe that for any $R \in (0,\infty)$
\bear\nonumber
t \int f^{2} \, (\widetilde{\log}_+ f)^{-2} 
&\le& 
t \int_{f \le R} f^{2} \, (\widetilde{\log}_+ f)^{-2}  + t \int_{f \ge R} f^{2} \,(\widetilde{\log}_+ f)^{-2} 
\\ \nonumber
&\le&
t \,  {R  \over(\widetilde{\log}_+ R)^{2} }
\int_{f \le R} f + { t  \over (\widetilde{\log}_+ R)^{2}} \int_{f \ge R} f^{2} 
\\  \label{eq:L2*to0}
&\le&
t \,  {MR  \over (\widetilde{\log}_+ R)^{2}}
 + { K \over (\widetilde{\log}_+ R)^{2}}   
\le
{M  + K \over (\widetilde{\log}_+ 1/t)^{2}}
 \to 0,
\eear
where we have used that $s \mapsto s/(\widetilde{\log}_+ s)^2$ is an increasing function in the second line, then 
the mass conservation and estimate \eqref{eq:L2tC} in the third line,  and we 
have chosen $R := t^{-1}$ in order to get the last inequality. We conclude to \eqref{eq:L43to0} by gathering 
\eqref{eq:L43bdd} and \eqref{eq:L2*to0}. 
\Black

\smallskip\noindent
{\sl Step 3.  Conclusion. } 
We consider two weak solutions $f_1$ and $f_2$ to the Keller-Segel equation \eqref{eq:KS} that we write in the mild form
$$
f_i(t) = e^{t\Delta} f_i(0) + \int_0^t e^{(t-s) \Delta } \nabla (V_i (s) \, f_i(s)) \, ds, 
\quad V_i = \KK * f_i,
$$
where $e^{t\Delta}$ stands for the heat semigroup defined in $\R^2$ by $e^{t\Delta} f := \gamma_t * f$, $\gamma_t(x) := (2\pi t)^{-1} \, \exp(-|x|^2/(2t))$. 
When we assume $f_1(0) = f_2(0)$, the difference $F := f_2 - f_1$ satisfies
$$
F(t) =  \int_0^t \nabla \cdot e^{(t-s) \Delta} (V_2(s) \, F(s)) \, ds
+ \int_0^t \nabla \cdot e^{(t-s) \Delta} (W(s) \, f_1(s)) \, ds = I_1 + I_2,
$$
with $W := V_2 - V_1$. For any $t > 0$, we define  
$$
Z_i(t) := \sup_{0 < s \le t} s^{1/4} \, \|  f_i(s) \|_{L^{4/3}}, 
\quad
 \Delta (t) := \sup_{0 < s \le t} s^{1/4} \, \|  F(s) \|_{L^{4/3}}.
$$
We then compute
\begin{eqnarray*}
J_1 
&:=& t^{1/4} \, \| I_1(t) \|_{L^{4/3}}
\\
&\le&
t^{1/4} \int_0^t \|  \nabla \cdot e^{(t-s) \Delta} (V_2 (s) \, F(s)) 
\|_{L^{4/3} } \, ds
\\
&\le& 
t^{1/4} \int_0^t  {C \over (t-s)^{3/4}} \, \| V_2 (s) \, F(s)  \|_{L^1 } \,
ds
\\
&\le& 
t^{1/4} \int_0^t  {C \over (t-s)^{3/4}} \, \| V_2 (s)  \|_{L^4}  \, \| F(s)
 \|_{L^{4/3} } \, ds
\\
&\le& 
t^{1/4} \int_0^t  {C \over (t-s)^{3/4}} \, \| f_2 (s)  \|_{L^{4/3}}  \, \|
F (s)  \|_{L^{4/3} } \, ds
\\
&\le& 
 \int_0^t  {C \over (t-s)^{3/4}} \, {t^{1/4}  \over s^{1/2}}  \, ds \, Z_2(t) \, \Delta(t)
 \\
&=& 
 \int_0^1  {C \over (1-u)^{3/4}} \, {du  \over u^{1/2}}   \, Z_2(t) \, \Delta(t), 
\end{eqnarray*}
where we have used the regularizing effect of the heat equation 
$$
\| \nabla (e^{t\Delta} g) \|_{L^{4/3}} \le \| \nabla \gamma_t \|_{L^{4/3}} \, \| g \|_{L^1} \le {C \over t^{3/4}} \, \| g \|_{L^1},
$$
at the third line, the H\"older inequality at the fourth line and the critical Hardy-Littlewood-Sobolev inequality \eqref{eq:HLSineqCritic} 
at the fifth line. 

Similarly, we have
\begin{eqnarray*}
J_2 
&:=& t^{1/4} \, \| I_2(t) \|_{L^{4/3}}
\\
&\le&
 \int_0^1  {C \over (1-u)^{3/4}} \, {du  \over u^{1/2}} \, \Delta(t)  \, Z_1(t) . 
\end{eqnarray*}
All together, we conclude thanks to \eqref{eq:L43to0} with the inequality
\bean
\Delta(t) 
&\le&   \int_0^1  {C \over (1-u)^{3/4}} \, {du  \over u^{1/2}}   \,  (Z_1(t) + Z_2(t))  \,  \Delta(t) \le {1 \over 2} \Delta(t)
\eean
for $t \in (0,T)$, $T > 0$ small enough, which in turn implies $\Delta(t) \equiv 0$ on $[0,T)$. 
\qed


 \section{Self-similar behaviour - Proof of Theorem~\ref{theo:LongTimeRate}} 
\label{sec:LongTimeBehave}
\setcounter{equation}{0}
\setcounter{theo}{0}

In this section we restrict ourself to the subcritical case $M < 8\pi$ and we investigate the self-similar long time
behaviour of generic solutions to the KS equation or more precisely, and equivalently, we investigate the long time
convergence to the self-similar profile of the rescaled solution $g$ defined through \eqref{eq:gt=}.
We start by recalling some known results on the self-similar profile and its stability. First,  we consider the
stationary problem \eqref{eq:statKSresc}.

\begin{theo}\label{theo:UniqGprofile} For any  $M  \in (0,8\pi)$, there exists a unique nonnegative self-similar profile $G = G_M$ of mass $M$ with finite second moment and finite entropy 
of the KS equation \eqref{eq:KS}, it is the unique solution to the stationary problem \eqref{eq:statKSresc} and it satisfies 
$$
G \in C^\infty(\R^2), \quad e^{-(1+\eps) |x|^2/2 + C_{1,\eps} } \le G \le e^{- (1-\eps) |x|^2/2 + C_{2,\eps} } ,
$$
for any $\eps \in (0,1)$ and some constants $C_{i,\eps} \in (0,\infty)$. 
Moreover, with the definitions \eqref{def:EE} of the modified free energy $\EE$ and \eqref{def:DDEE} of  the modified dissipation of the free energy $\DD_\EE$, the self-similar profile $G$ is characterized as the unique solution to the optimization problem
\beqn\label{eq:minimizeEE}
\tilde g \in \ZZ_M, \quad \EE(\tilde g) = \min_{g \in \ZZ_M} \EE(g), 
\eeqn
where $\ZZ_M := \{ g \in L^1_+ \cap L^1_2; \,\, M_0(g) = M \}$, as well as the unique function $g \in \ZZ_M$ such that $\DD_\EE(g) = 0$. 
\end{theo}

That theorem follows by a combination of known results. On the one hand, as a consequence of the fact that $U := - \KK * G$ satisfies  \eqref{eq:U=exp}
together with  the elementary inequality 
\beqn\label{eq:borneU}
\forall \, x \in \R^2 \quad \bigl| U(x) + { M \over 2\pi} \, (\log |x|)_+ \bigr| \le C,
\eeqn
where $C$ only depends on $M$, $M_2(G)$ and $\HH(G)$ (see \cite[Lemma 23]{MR2226917} and the argument presented in order
to bound $\RR_\eps(g)$ in step 2 of the proof of Theorem~\ref{theo:ApostUnifSm}), and the Naito's variant  \cite{MR1972310} 
of the famous Gidas, Ni, Nirenberg radial symmetry result on solutions to Poisson type equations, it has been established in 
\cite[Lemma 25]{MR2226917} that  $U$ is radially symmetric. It follows that any self-similar profile $G$ is radially symmetric. 
On the other hand, the uniqueness of  radially symmetric  self-similar profiles has been proved in   \cite[Theorem 3.1]{BKLN} 
(see also \cite[Theorem 1.2]{MR2929020}) and that concludes the proof of the uniqueness of the solution  to the stationary problem \eqref{eq:statKSresc}.
The smoothness property 
is established in \cite[Lemma 25]{MR2226917}  and the behaviour for large values of $|x|$ is a immediate consequence of \eqref{eq:borneU}.
It is clear from \eqref{eq:dissipRescE} that any solution $\tilde g$ to the minimization problem \eqref{eq:minimizeEE} also satisfies $\DD_\EE(\tilde g) = 0$
which in turns implies that $\log \tilde g+ |x|^2/2 + \kappa * \tilde g = 0$ and then $\tilde g$ is a solution to the  stationary problem \eqref{eq:statKSresc}. 

\smallskip

Second, the profile $G$ is a stationary solution to the evolution equation  \eqref{eq:KSresc} and 
the associated linearized equation reads 
 $$
  \partial_t h = \Lambda h := \hbox{div}_x \bigl( \nabla h + x \,h + (\KK * G) \, h + (\KK*h) \, G). 
$$
We briefly explain the spectral analysis of $\Lambda$ in the Hilbert space $E :=L^2(G^{-1/2})$ {\it of self-adjointness}
performed in \cite{CamposDolbeault2012}. 
{
Defining $h_{0,0} := {\partial  G_M /\partial M}$, it is (formally) clear that  $h_{0,0}$ is a first eigenfunction of the operator $\Lambda$
 associated to the first eigenvalue $\lambda = 0$,  and it has  been furthermore shown in \cite[Lemma 8]{CamposDolbeault2012} that 
 the null space $N(\Lambda) = \hbox{vect} (h_{0,0})$. Moreover, defining the bilinear form 
$$
\langle f,g \rangle := \int_{\R^2} f \, g  \, G^{-1} \, dx +  \int_{\R^2} \!  \int_{\R^2} f(x) \, g(y) \, \kappa (x-y) \, dx dy,
$$
and the associated quadratic form $Q_1[f] := \langle f, f \rangle$, it has been shown in \cite[Section 4.3]{CamposDolbeault2012} that $Q_1$ is nonnegative, that $Q_1[h_{0,0}] = 0$ and that 
$$
Q_1[f] = 0 \ \hbox{and}\  \langle f,h_{0,0} \rangle = 0
\quad\hbox{imply}\quad f = 0.
$$
%
As a consequence $Q_1[\cdot]$ defines an Hilbert norm on the  linear submanifold 
$$
E_0^\perp :=   \{ f \in E; \,\,   \langle   f,h_{0,0} \rangle = 0 \} =   \{ f \in E; \,\,  M( f ) = 0 \} 
$$
 which is equivalent to the initial norm $\| \cdot \|_E$. That new norm is suitable for exhibiting a spectral gap for the operator $\Lambda$ and to make the stability analysis of the associated semigroup $e^{t\Lambda}$. }

\begin{theo}[\cite{CamposDolbeault2012}] \label{theo:SGL2}For any $g \in E_0^\perp$ which belongs to the domain of $\Lambda$,  there holds
\beqn\label{eq:Q2Q1}
\langle \Lambda g,g\rangle  \le - Q_1[ g ].
\eeqn
Moreover, there exists $a^* < -1$ and $C > 0$ so that 
\beqn\label{eq:SgE}
\| e^{t\Lambda} h  - e^{-t} \, \Pi_1 h - \Pi_0 h \|_E \le C \, e^{a^*t} \, \| h - (\Pi_1 + \Pi_0) h \|_E
\qquad \forall \, t \ge 0, \,\, \forall \, h \in E, 
\eeqn
where $\Pi_0$ is the ($Q_1$-orthogonal) projection on $\hbox{Vect}(h_{0,0})$, 
also defined as $\Pi_0 h := M(h) \, h_{0,0}$, \Black and $\Pi_1$ is the 
($Q_1$-orthogonal) projection on $\hbox{Vect}(h_{1,1},h_{1,2})$ where $h_{1,i} := \partial_{x_i} G$. 
\end{theo}

Inequality \eqref{eq:Q2Q1} is nothing but \cite[Theorem 15]{CamposDolbeault2012} and \eqref{eq:SgE} is 
a consequence of the fact that the spectrum of $\Lambda$ is discrete and included in the real line  and that the second (larger)
eigenvalue of $\Lambda$ is $-1$, see \cite[Section 4]{CamposDolbeault2012}.

\smallskip
Our first main result in this section is a linearized stability result in a large space $\EE$, namely we consider
\bean
\EE &:=& L^{4/3}_k(\R^2), \quad k > 3/2.
\eean
We consider that space because it is the larger space in terms of moment decay in which we are able to 
prove a (optimal) spectral gap on the linearized semigroup. For such a general Banach space framework and the
associated spectral analysis issue, we adopt the classical notations of \cite{Pazy,Kato} used in \cite{GMM}, for more details we refer to \cite[Section 2.1]{GMM}
and the references therein (in particular \cite{Kato,Pazy,MR1721989}).

\begin{theo}\label{theo:SGL43} 
For any $k > 3/2$ and any $a > \bar a := \max(a^*,a(k))$, $a(k) := 1/2 - k$ (so that $a(k) < -1$) there exists a constant $C_{k,a}$ so that 
$$
\| e^{t\Lambda} h - e^{-t} \, \Pi_1 h - \Pi_0 h \|_\EE \le C \, e^{at} \, \| h - \Pi_1 h - \Pi_0 h  \|_\EE
\quad \forall \, t \ge 0, \,\, \forall \, h \in \EE,
$$
 where again $\Pi_0$ stands for projection on the eigenspace $\hbox{Vect}(h_{0,0})$ associated to the eigenvalue $0$
 and $\Pi_1$ stands for projection on the eigenspace  $\hbox{Vect}(h_{1,1},h_{1,2})$ associated to the eigenvalue $-1$. 
 Both operators are defined through the Dunford formula (see \cite[Section 2.1]{GMM} or better \cite[III-(6.19)]{Kato})
 $$
 \Pi_{\xi} := - {1 \over 2i\pi} \int_{ |z - \xi| = r  } (\Lambda - z)^{-1} \,  dz, \quad \xi = 0, -1, \,\,\, r > 0 \, \, (\hbox{small enough}), 
 $$
 but also in a simpler manner $\Pi_0 h = M(h) \, h_{0,0}$ for any $h \in \EE$. 
\end{theo}

The proof is a straightforward adaptation of arguments  of {\it ``functional extension of semigroup spectral gap estimates" } developed in \cite{GMM} for
the Fokker-Planck equation.  

\begin{lem}\label{lem:dissipB} For any $k \ge 0$ fixed,  there exists a constant $C_k$ such that for any $g \in D(\Lambda)$, there holds
\beqn\label{eq:dissipB}
\langle \Lambda g, g^\dagger \rangle_\EE \le C_k \int |g|^{4/3} \, \langle x \rangle^{{\frac 43} k-1} + 
 \bigl( \frac12 - k  \bigr) \int |g|^{4/3} \, \langle x \rangle^{{\frac 43}k} ,
\eeqn
where $g^\dagger := \bar g \, |g|^{-2/3}$ (here $\bar g$ stands for the complex conjugate of $g$). 

\end{lem}

\noindent
{\sl Proof of Lemma~\ref{lem:dissipB}. } For the sake of simplicity we assume $g \ge 0$ so that $g^\dagger = g^{1/3}$, we set $\ell := 4k/3$, we write 
$$
\langle \Lambda g, g^\dagger \rangle_\EE = \int_{\R^2} (\Lambda g) \, g^{1/3} \, \langle x \rangle^\ell = T_1 + ... + T_4,
$$ 
and we compute each term $T_i$ separately. First, performing two integrations by part, we have 
\bean
T_1 
&:=& \int_{\R^2} (\Delta g) \, g^{1/3} \, \langle x \rangle^\ell \, dx
\\
&=& -{1 \over 3} \int_{\R^2} |\nabla g |^2 \, g^{-2/3}  \, \langle x \rangle^\ell \, dx + \frac34\int_{\R^2}   g^{4/3}  \, \Delta \langle x \rangle^\ell \, dx .
\eean 
 Second, performing one integration by part, we have 
\bean
T_2 
&:=& \int_{\R^2} ( 2 g + x \cdot \nabla g) \, g^{1/3} \, \langle x \rangle^\ell \, dx
\\
&=&  \int_{\R^2} \bigl\{ \frac12 \langle x \rangle^{\ell-2} + \bigl( \frac12 - k \bigr) \langle x \rangle^{\ell} \bigr\} \, g^{4/3} \, dx.
\eean 
Third, performing one integration by part, we have 
\bean
T_3 
&:=& \int_{\R^2} \bigl(  2 \, G g + (\KK * G) \cdot \nabla g \bigr)   \, g^{1/3} \, \langle x \rangle^\ell \, dx
\\
&=& \frac54 \int_{\R^2}  G \, g^{4/3} \, \langle x \rangle^\ell \, dx
- \frac34 \int_{\R^2} g^{4/3} \,  (\KK * G) \cdot \nabla_x   \langle x \rangle^\ell \, dx
\\
&\le&C  \int_{\R^2}   g^{4/3} \,     \langle x \rangle^{\ell -1}\, dx,
\eean
for some constant $C \in (0,\infty)$. 

\noindent
Fourth and last, thanks to the H\"older inequality and  the critical Hardy-Littlewood-Sobolev inequality \eqref{eq:HLSineqCritic}, 
we have
\bean
T_4
&:=& \int_{\R^2}   (\KK * g) \cdot \nabla G   \,  g^{1/3} \, \langle x \rangle^\ell \, dx
\\
&\le&  \| \nabla G \, \langle x \rangle^k \|_\infty \, \| g \|_{L^{4/3}}^{1/3} \, \| \KK * g \|_{L^4} 
\le C\, \| g \|_{L^{4/3}}^{4/3}. 
\eean
Gathering all these estimates, we get \eqref{eq:dissipB}. \qed

 \medskip
 We define 
 $$
 \AA g := N \chi_R \, g \quad \hbox{and} \quad \BB g = \Lambda g - \AA \, g,
 $$
 for some truncation function $\chi_R(x) :=\chi(x/R)$, $\chi \in \DD(\R^2)$, ${\bf 1}_{B(0,1)} \le \chi \le {\bf 1}_{B(0,2)}$, and
 some constants $N,R >0$. 
 
\smallskip  
 We clearly have 
 \beqn\label{eq:boundAA}
 \AA \in \BBB(L^2,E) \subset \BBB(E) \quad \hbox{ and } \quad \AA \in \BBB(L^{4/3},\EE) \subset \BBB(\EE).
 \eeqn
 
\smallskip
From lemma~\ref{lem:dissipB} we easily have that for any $a > a(k)$ there exist $N$ and $R$ large enough so that 
 $\BB-a$ is dissipative in $\EE$ (see \cite[Chapter I, Definition 4.1]{Pazy})  in the sense that 
 \beqn\label{eq:boundBB}
 \langle g^*, (\BB -a) \, g \rangle_{\EE',\EE}   \le 0,
\eeqn
 where $g^* := \bar g \, |g|^{-2/3} \, \| g \|_\EE^{2/3} \in \EE'$. 
 
%
%

\Black
 \begin{lem}\label{lem:Nash} There exist some constants $C > 0$ and $b \in \R$ such that the semigroup $S_\BB(t) = e^{\BB t}$ satisfies
\beqn\label{eq:NashB}
 \| S_\BB(t) h \|_{L^2_1} \le {C \, e^{bt} \over t^{1/2}} \, \| h \|_{L^{4/3}_1} \quad \forall \, h \in L^{4/3}_1, \,\,\forall \, t > 0.
\eeqn
 \end{lem}
 
 \noindent
{\sl Proof of Lemma~\ref{lem:Nash}. } The proof of the hypercontractivity property as stated in Lemma~\ref{lem:Nash}
is a classical consequence of the Gagliardo-Nirenberg  inequality. For the sake of completeness we sketch it. Arguing similarly as in the proof
of Lemma~\ref{lem:BdLp} and Lemma~\ref{lem:dissipB} and denoting $h_t := e^{t\BB} h$, we compute 
%
%
\bean
&&\frac12 {d \over dt} \int |h_t|^2 \langle x \rangle^2 
= - \int |\nabla h_t \langle x \rangle |^2
+ \int h_t \, (\KK*h_t) \cdot \nabla G  \langle x \rangle^2
\\
&&\qquad+ \int |h_t|^2 \, \Bigl\{ 1 -   |\nabla  \langle x \rangle |^2  +  \langle x \rangle^2  \, \bigl(\frac 32 G - N \, \chi_R \bigr) 
+ \langle x \rangle \, \KK*G \cdot \nabla \langle x \rangle \Bigr\}. 
\eean
On the one hand, 
thanks to the Gagliardo-Nirenberg inequality \eqref{eq:GNineq} with $q=4/3$, $r=2$ and $a=1/3$, we know that 
 $$
 \int |\nabla (h \langle x \rangle) |^2 \ge C_{GN}^{-6} \, 
 \Bigl( \int | h \, \langle x \rangle |^2 \Bigr)^3 \,  \Bigl( \int | h \, \langle x \rangle |^{4/3} \Bigr)^{-3}. 
 $$
On the other hand, introducing the splitting $\KK = \KK_0 + \KK_\infty$ with $\KK_0 := \KK \, {\bf 1}_{|z| \le 1}$ and $\KK_\infty := \KK \, {\bf 1}_{|z| \ge 1}$ and 
using the H\"older inequality and the Young inequality, we have 
\bean
 \int h_t \, (\KK*h_t) \cdot \nabla G \langle x \rangle^2
 &\le& \| \nabla G  \, \langle x \rangle^2 \|_{L^\infty} \, \| h \|_2 \, \| \KK_0 * h \|_{L^2} +\| \nabla G  \, \langle x \rangle^2 \|_{L^2} \, \| h \|_{L^2} \, \| \KK_\infty * h \|_{L^\infty} 
\\
&\le& C \, ( \| \KK_0 \|_{L^1} \,  \| h \|_{L^2}^2 + \| h \|_{L^2}  \,  \| \KK_0 \|_{L^3} \|  h \|_{L^{3/2}} )  
 \\
&\le& C \,  \| h \|_{L^2_1}^2  . 
 \eean
 We also bound the last term by $C \, \| h_t \|_{L^2_1}^2$. 
All together and using the notations $X(t) := \| h_t \|_{L^2_1}^2$ and $Y(t) := \| h_t \|_{L^{4/3}_1}^{4/3}$ and the fact that $Y(t) \le Y(0)$ thanks to Lemma~\ref{lem:dissipB}, we get
$$
X' \le  - \alpha \, (X/Y(0))^3 + \beta \, X
$$
for some constants $\alpha, \beta > 0$. The estimate \eqref{eq:NashB} is then a classical consequence to the above differential inequality. 
\qed

\medskip \noindent
{\sl Proof of Theorem~\ref{theo:SGL43}. } We immediately deduce from \eqref{eq:boundAA} and Lemma~\ref{lem:Nash} that 
$$
\| \AA S_\BB(t) \|_{\BBB(\EE,E)} \le {C \over t^{1/2}} \, e^{bt} \quad \forall \, t > 0, 
$$
for some constants $C > 0$ and $b \in \R$. As a consequence, proceeding as in   \cite[section 3]{GMM} or \cite[Lemma~2.4]{MM*}, we deduce that the time convolution function $(\AA S_\BB)^{(*\ell)}$ defined iteratively by $(\AA S_\BB)^{(*1)} := (\AA S_\BB)$, $(\AA S_\BB)^{(*\ell)} := (\AA S_\BB)^{(*(\ell-1))}  * (\AA S_\BB)$, for any $\ell \ge 2$,  satisfies 
\beqn\label{eq:ASB2}
\| (\AA S_\BB)^{(*\ell)} (t) \|_{\BBB(\EE,E)} \le C_\ell \, e^{b_\ell t} \quad \forall \, t > 0, 
\eeqn
for some constants $C_\ell> 0$ and $b_\ell \in (a(k),-1)$ for $k > 3/2$ and $\ell \ge 2$ large enough. Putting together Theorem~\ref{theo:SGL2} and the properties \eqref{eq:boundAA}, \eqref{eq:boundBB}
and \eqref{eq:ASB2} we observe that $\Lambda = \AA + \BB$ satisfies all the assumptions of \cite[Theorem~2.13]{GMM}. As a consequence, the conclusions of Theorem~\ref{theo:SGL2} hold true by a straightforward application of \cite[Theorem~2.13]{GMM}. 
\qed

 \smallskip
 
 Before going to the proof of Theorem~\ref{theo:LongTimeRate}  we present two results that will be useful during the proof of that Theorem. 
 
 \begin{lem}\label{lem:FDF} For any $M \in (0,8\pi)$, $k' >2 \ge k >3/2$, $M_{k'} \ge (k'-1)^{k'/2} \, M$ and $\CC > 0$,  there exists an increasing function $\eta : [0,\infty) \to [0,\infty)$,
 $\eta(0) = 0$, $\eta(u) > 0$ for any $u > 0$, such  that 
\beqn\label{eq:FDF}
 \forall \, g \in \ZZ \qquad \DD_\EE(g) \ge \eta( \| g - G \|_{L^{4/3}_k} )
\eeqn
 where 
  $$
 \ZZ := \{g \in L^1_+ (\R^2), \,\, M(g) = M, \,\, M_{k'}(g) \le M_{k'}, \,\, \| g \|_{W^{2,\infty}} \le \CC \}.
 $$
 
 \end{lem}
 
%
%
 
\noindent
{\sl Proof of Lemma~\ref{lem:FDF}. }  We proceed by contradiction. If \eqref{eq:FDF} does not hold, there exists 
a sequence $(g_n)$ in $\ZZ$ and a real $\delta > 0$ such that 
$$
\DD_\EE(g_n) \to 0 \,\,\hbox{as}\,\, n \to 0
\quad\hbox{and}\quad
 \| g - G \|_{L^{4/3}_{k}}  \ge \delta.
$$
Therefore, on the one hand, there exists $\bar g \in \ZZ$ such that, up to the extraction of the subsequence, there holds
$g_n \to \bar g$ strongly in $L^{4/3}_{k}$, so that $\| \bar g - G\|_{L^{4/3}_k} \ge \delta$. 
Using again  $g_n \to \bar g$
and the critical Hardy-Littlewood-Sobolev inequality \eqref{eq:HLSineqCritic}, we deduce that $\sqrt{g_n} \, \KK * g_n \to  \sqrt{\bar g} \, \KK * \bar g $
strongly in $L^1_{loc}(\R^2)$ and then $2\nabla \sqrt{g_n} + \sqrt{g_n} \, \KK * g_n \wto 2\nabla \sqrt{\bar g} + \sqrt{\bar g} \, \KK * \bar g $  in $\DD'(\R^2)$.
Since $(\nabla \sqrt{g_n} + \sqrt{g_n} \, \KK * g_n)$ is 
bounded in $L^2$, that implies that \Black
$2\nabla \sqrt{g_n} + \sqrt{g_n} \, \KK * g_n \wto 2\nabla \sqrt{\bar g} + \sqrt{\bar g} \, \KK * \bar g $ weakly in $L^2(\R^2)$
and then 
$$
\DD_\EE(\bar g) = \| 2\nabla \sqrt{\bar g} + \sqrt{\bar g} \, x + \sqrt{\bar g} \, \KK * \bar g\|_{L^2}^2 \le \liminf \DD_\EE(g_n) = 0.
$$
We easily conclude thanks to the mass condition $M_0(\bar g) = M$  and the uniqueness Theorem~\ref{theo:UniqGprofile} that $\bar g = G$. That is our contradiction.
%
  \qed

\begin{lem}\label{lem:LambdaDissipatif} 
Define $\EE_2 := R(I-\Pi_0 - \Pi_1)$ the supplementary linear submanifold to the eigenspaces associated to the 
eigenvalues $0$ and $-1$. There exists a norm $\Nt \cdot \Nt$ on $\EE_2$ equivalent to the initial one
$\| \cdot \|_\EE$  so that 
\beqn\label{ineq:NormeN*}
{d \over dt} \Nt e^{t\Lambda}  f \Nt^2 \le 
- 2 \,   \Nt  e^{t\Lambda} f  \Nt ^2  \quad \forall \, t \ge 0, \,\, \forall \, f \in \EE_2.
\eeqn

  \end{lem}

\noindent
{\sl Proof of Lemma~\ref{lem:LambdaDissipatif}. } This result is nothing but   \cite[Proposition 5.14]{GMM}. For the sake of 
completeness and because we will need to use the same computations at the nonlinear level, we just check it below. 
First recall that from  Theorem~\ref{theo:SGL43}, we know that 
for any $a \in ( \bar a,-1)$ there exists $C = C(a)$ such that 
$$
\| e^{\Lambda t} \, f \|_\EE \le C \, e^{a \, t} \, \| f \|_\EE, \quad \forall \, t \ge 0,  \,\, \forall \, f \in \EE_2, 
$$
and on the other hand, from Lemma~\ref{lem:dissipB} there exists some constant 
$b \in \R$   such that 
$$
\langle \Lambda f, f^* \rangle \le b \, \| f \|_\EE^2.
$$
We define 
\beqn\label{eq:defNormNt}
\Nt f \Nt^2  := \eta \, \| f \|^2_\EE + \int_0^\infty \| e^{\tau \Lambda} \, e^{ \tau} \, f \|^2_\EE  \, d\tau 
\eeqn
with $\eta \in (0, (b+1)^{-1})$. The norm  $\Nt \cdot \Nt$ is clearly well defined and it is equivalent to $\| \cdot \|_\EE$ because 
$$
\forall \, f \in \EE_2, \quad
\eta \, \| f \|^2_\EE  \le   
\eta \, \| f \|_\EE^2 + \int_0^\infty \| e^{\Lambda \tau}  \, e^{ \tau} \, f \|_\EE^2  \, d\tau
\le \Bigl( \eta \,  + \int_0^\infty  C^2 \, e^{2 \, (a+1) \tau} \, d\tau \Bigr) \, \|  f \|_\EE^2 . 
$$
Next, for $f \in \EE_2$ and with the notation $f_t := e^{\Lambda t} f$, we compute
\bean
{d \over dt} \Nt e^{\Lambda t} f \Nt^2  
&=&
 \eta \, {d \over dt} \| f_t \|^2 + \int_0^\infty  {d \over dt} \| e^{\Lambda \, (t+\tau) + \tau} \, f \|^2  \, d\tau
\\
&=&
2 \eta \, \langle f^*_t, \Lambda f_t \rangle 
+ \int_0^\infty  
\Bigl\{ {d \over d\tau} \| e^{\Lambda \, (t+\tau) + \tau } \, f \|^2  - 2  \,  \| e^{\Lambda \, (t+\tau) + \tau } \, f \|^2   \Bigr\} \, d\tau
\\
&\le&
2 \eta \, b \, \| f_t \|^2 +  \Bigl[  \| e^{\Lambda \, (t+\tau) + \tau } \, f \|^2 \Bigr]_0^\infty  - 2  \int_0^\infty   \| e^{\Lambda \tau}\, e^{\tau } \, f_t \|^2 \, d\tau
\\
&=&
\Bigl\{ 2 \eta \, (b-a) - 1 \Bigr\}  \| f_t \|^2 - 2 \, \Bigl\{ \eta \| f_t \|^2 +  \int_0^\infty   \| e^{\Lambda \tau} \, e^{ \tau } \, f_t \|^2 \, d\tau \Bigr\}
\\
&\le& - 2 \, \Nt e^{\Lambda t} f \Nt^2 ,
\eean
so that \eqref{ineq:NormeN*} is proved.  \qed

\medskip
We conclude with the proof of the long time convergence result. 

\medskip\noindent
{\sl Proof of Theorem~\ref{theo:LongTimeRate}. } The proof follows the same strategy as in \cite{Mcmp,MMcmp,GMM} 
(see also \cite{MR900501,MR946973,MR1338453} where similar proof is carried on in the context of the Boltzmann equation). 
We split the proof into four steps. 

\medskip\noindent
{\sl Step 1. }  We consider a solution $g$ to the rescaled equation \eqref{eq:KSresc} with initial datum $f_0 \not=G$.
Thanks to Theorem~\ref{theo:ApostUnifSm} there holds $g(t) \in \ZZ$ for any $t \ge 1$. 
For any $\delta > 0$ and $T := (\EE(f_0) - \EE(G))/ \eta^{-1}(\delta)+1$ there exists $t_0 \in [1,T]$
so that 
\beqn\label{eq:DEg<d}
\DD_\EE(g(t_0)) \le \eta^{-1}(\delta)
\eeqn
because on the contrary we would have from \eqref{eq:dissipRescE}
$$
{d \over dt} (\EE(g(t)) - \EE(G)) \le -  \eta^{-1}(\delta) \quad \hbox{on} \quad (1,T), 
$$
and then 
$$
\EE(g(T)) - \EE(G) \le - (\EE(f_0) - \EE(G))  < 0
$$
which is in contradiction with the fact that $G$ satisfies
$
\EE (G) < \EE(f)$ $\forall \, f \in \ZZ \backslash \{ G \}$ from Theorem~\ref{theo:UniqGprofile}. 
We deduce from \eqref{eq:DEg<d} and Lemma~\ref{lem:FDF} that 
$$
\| g(t_0) - G \|_{L^{4/3}_{k'}} \le \delta.
$$
 
\medskip\noindent
{\sl Step 2. }  The function $h := g - G$ satisfies the equation
$$
\partial_t h = \Lambda h + \hbox{div} ( h \, \KK * h).
$$
We introduce the splitting 
$$
h = h_0 + h_1 + h_2, \quad h_{12} = h_1 + h_2
$$
with 
$$
h_0 := \Pi_0 h, 
\quad
h_1 := \Pi_1 h ,
$$
so that the evolution of $h_1$ and $h_2$ are given by 
\beqn\label{eq:h1}
\partial_t h_1 = - h_1 + \Pi_1 [ \hbox{div} ( h \, \KK * h) ] 
\eeqn
and
\beqn\label{eq:h2}
\partial_t h_2 = \Lambda h_2 + \QQ_2, \quad \QQ_2 := \Pi_2[ \hbox{div} ( h \, \KK * h) ] . 
\eeqn
Because of the mass conservation $M(g(t)) = M(G)$, there holds $h_0(t) = \Pi_0 h(t) = h_{0,0} \, M(h(t)) = 0$.  \Black
Moreover, from \eqref{eq:h1} and with the notation $h^*_1 = h_1 \, |h_1|^{-1/3} \, \| h_1 \|_{L^{4/3}_k} ^{2/3}$, we clearly have 
\bear\nonumber
{d \over dt} \| h_1 \|_{L^{4/3}_k} ^2  
&=& 2 \,  \langle- h_1 + \Pi_1 [ \hbox{div} ( h \, \KK * h) ] , h_1^* \rangle
\\ \nonumber
&\le& -2 \,  \| h_1  \|^2_{L^{4/3}_k} + 2 \,  \| h_1  \|_{L^{4/3}_k} \,  \|  \Pi_1 [ \hbox{div} ( h \, \KK * h) ] \|_{L^{4/3}_k}
\\ \label{eq:estimh1}
&=& - 2 \, \| h_1  \|_{L^{4/3}_k}^2 + C \,\| h_1  \|_{L^{4/3}_k} \,  \|  \hbox{div} ( h \, \KK * h)  \|_{L^{4/3}_k}. 
\eear
 
\medskip\noindent
{\sl Step 3. Estimate on the nonlinear term. } We make the splitting 
$$
\| \hbox{div} (h \, \KK * h ) \|_{L^{4/3}_k} \le I_1 + I_2, 
\quad I_1 := \| h^2 \|_{L^{4/3}_k}, 
\quad I_2 := \|  \nabla h \cdot \KK * h\|_{L^{4/3}_k}, 
$$
and we compute each term separately. On the one hand, using the H\"older inequality and 
the  Galgliardo-Nirenberg inequality (see \cite[Chapter IX, inequality (86)]{BrezisBook}) in
dimension 2
$$
\| u \|_{L^p} \le C \, \| u \|_{L^q}^{1-a} \, \| u \|_{W^{1,r}}^a, 
\quad a = {{1 \over q} - {1 \over p} \over {1 \over q} + {1 \over 2}-{1\over r}} , 
$$
with $r=p=\infty$, $q=4/3$ and $a= 3/5$, 
we have 
\bean
I_1 
&\le& \| h \|_{L^\infty}   \, \| h \|_{L^{4/3}_{k}}
\le  C \, \|  h \|_{L^{4/3}_k}^{7/5} \, \| h \|_{W^{1,\infty}}^{3/5}. 
\eean
On the other hand, thanks to  the critical Hardy-Littlewood-Sobolev inequality \eqref{eq:HLSineqCritic}, 
 the elementary inequality 
$$
\| \nabla u \|_{L^2_k}^2 = - \int_{\R^2} u \, \hbox{div} ( \langle x \rangle^{2k} \, \nabla u ) \le  C \, \| u \|_{W^{2,\infty}} \, \| u \|_{L^1_{2k}},
$$
and the H\"older inequality 
$$
 \| u \|_{L^1_{2k}} \le  \| \langle x \rangle^{-1} \|_{L^{4\gamma/\alpha}}^{1/\gamma} \, \| u \|_{L^{4/3}_{k}} ^{\alpha} \, \| u \|_{L^1_{k'}} ^{1-\alpha} , 
 $$
 with $0<\alpha<1$, $2\gamma > \alpha$ and $k' = k'(\alpha,\gamma) := ((2-\alpha)k + \gamma)/(1-\alpha)$,  
we have 
\bean
I_2 
&\le& \| \nabla h \|_{L^{2}_k}\, \| \KK * h  \|_{L^4}
\le  C_{\alpha,\gamma} \, \|  h \|_{L^{4/3}_k}^{1+\alpha/2} \, \| h \|_{W^{2,\infty}}^{1/2} \, \| h \|^{(1-\alpha)/2}_{L^1_{k'}}.
\eean
To make the computtaions simpler, when  $k' = 4$, we can take $k=8/5 > 3/2 $, $\gamma/\alpha = 5/8 > 1/2$ and we get $\alpha = 32/121\in (0,1)$ and $\alpha/2 < 2/5$. 
All together we find
$$
\forall \, h \in \ZZ, \qquad \| \hbox{div} (h \, \KK * h ) \|_{L^{4/3}_k} \le C \, \| h \|_{L^{4/3}_k}^{1+\alpha/2}. 
$$
Thanks  to Theorem~\ref{theo:ApostUnifSm} we have $h(t) \in \ZZ$ for all $t \ge 1$ (where in the definition $\ZZ$ the constant $\CC$ is given by \eqref{eq:bddW2infty}), and we conclude with 
\beqn\label{eq:NLestim}
\forall \, t \ge 1, \qquad \| \hbox{div} (h(t) \, \KK * h(t) ) \|_{L^{4/3}_k} \le C \, \| h(t) \|_{L^{4/3}_k}^{1+\alpha/2}.
\eeqn
 It is worth noticing that in the limit  $k \to 3/2$, $\gamma/\alpha \to 1/2$ and $\alpha \to 0$, we find $k'= 3$. In other words, one can easily verify that \eqref{eq:NLestim} still holds for any $k' > 3$ (with another choice of $\alpha \in (0,1)$). 
 
\medskip\noindent
{\sl Step 4. Estimate on the remaining term and conclusion. } 
From \eqref{eq:h2}, using the norm $\Nt \cdot \Nt$ defined in \eqref{eq:defNormNt} and the notation $S_\tau := e^{\tau \Lambda} \, e^\tau$, 
we compute 
\bear\nonumber
{d \over dt}  \Nt h_2 \Nt^2
&=& 
 \eta \, \langle h_2^*,\Lambda h_2 \rangle + \int_0^\infty \langle (S_\tau h_2)^*, S_\tau \Lambda h_2 \rangle   \, d\tau 
\\ \nonumber
&&+
 \eta \, \langle h_2^*, \QQ_2 \rangle + \int_0^\infty \langle (S_\tau h_2)^*, \QQ_2 \rangle   \, d\tau 
\\ \label{eq:estimh2}
&\le& - 2 \, \Nt h_2 \Nt^2 + C \, \| h_2 \|_{L^{4/2}_k} \, \| \hbox{div} (h \, \KK * h ) \|_{L^{4/2}_k}, 
\eear
where we have used Lemma~\ref{lem:LambdaDissipatif} in order to bound the first (linear) term and the equivalence between the 
two norms $\Nt \cdot \Nt$ and $\| \cdot \|_{L^{4/2}_k}$ in order to estimate the second one (which involves the nonlinear quantity). 
Gathering \eqref{eq:estimh1}, \eqref{eq:estimh2},  \eqref{eq:NLestim}, we clearly see that 
$$
u(t) :=  \| h_1 \|_\EE^2+ \Nt h_ 2 \Nt^2 
$$
satisfies the differential inequality
$$
u' \le -2u + C \, \| h \|^{2+\alpha}  \quad \hbox{on}\,\, (0,\infty) , 
$$
and then thanks to 
the first step
\beqn\label{eq:dudt1}
u' \le -2u + C \, u^{1+\alpha/2}  \quad \hbox{on}\,\, (t_0,\infty), \quad u(t_0) \le K_2 \, \delta. 
\eeqn
Taking $\delta > 0$ small enough in the first step, we classically deduce that 
\beqn\label{eq:dudt2}
u(t)  \le C_{a} \, e^{2a \, t }  \quad \forall \, t \ge t_0
\eeqn
for any $a > -1 $, so that for $a$ close enough to $-1$, we deduce from \eqref{eq:dudt1}-\eqref{eq:dudt2} that
$$
u' \le -2u + K_2 \, e^{-2t}  \quad \hbox{on}\,\, (t_0,\infty),  
$$
from which we easily conclude $u \le C \, e^{-2t}$. 
 \qed

\bibliographystyle{acm}


\signge \signsm  

\end{document}